\pdfoutput=1
\RequirePackage{ifpdf}
\ifpdf % We~are running pdfTeX in pdf mode
\documentclass[pdftex]{sigma}
\else
\documentclass{sigma}
\fi

\numberwithin{equation}{section}

\newtheorem{Theorem}{Theorem}[section]
\newtheorem{Corollary}[Theorem]{Corollary}
\newtheorem{Lemma}[Theorem]{Lemma}
\newtheorem{Proposition}[Theorem]{Proposition}
{ \theoremstyle{definition}
\newtheorem{Definition}[Theorem]{Definition}}

\usepackage{mathtools}

\begin{document}
\allowdisplaybreaks

\newcommand{\arXivNumber}{2204.09206}

\renewcommand{\PaperNumber}{079}

\FirstPageHeading

\ShortArticleName{Noncolliding Macdonald Walks with an Absorbing Wall}

\ArticleName{Noncolliding Macdonald Walks\\ with an Absorbing Wall}

\Author{Leonid PETROV}

\AuthorNameForHeading{L.~Petrov}

\Address{University of Virginia, Charlottesville, VA, USA}
\Email{\href{mailto:lenia.petrov@gmail.com}{lenia.petrov@gmail.com}}
\URLaddress{\url{https://lpetrov.cc/}}

\ArticleDates{Received June 07, 2022, in final form October 16, 2022; Published online October 20, 2022}

\Abstract{The branching rule is one of the most fundamental properties of the Macdonald symmetric polynomials. It expresses a Macdonald polynomial as a nonnegative linear combination of Macdonald polynomials with smaller number of variables. Taking a limit of the branching rule under the principal specialization when the number of variables goes to infinity, we obtain a Markov chain of $m$ noncolliding particles with negative drift and an absorbing wall at zero. The chain depends on the Macdonald parameters $(q,t)$ and may be viewed as a discrete deformation of the Dyson Brownian motion. The trajectory of the Markov chain is equivalent to a certain Gibbs ensemble of plane partitions with an arbitrary cascade front wall. In the Jack limit $t=q^{\beta/2}\to1$ the absorbing wall disappears, and the Macdonald noncolliding walks turn into the $\beta$-noncolliding random walks studied by Huang~[\textit{Int. Math. Res. Not.} \textbf{2021} (2021), 5898--5942, arXiv:1708.07115]. Taking $q=0$ (Hall--Littlewood degeneration) and further sending $t\to 1$, we obtain a~continuous time particle system on $\mathbb{Z}_{\ge0}$ with inhomogeneous jump rates and absorbing wall at zero.}

\Keywords{Macdonald polynomials; branching rule; noncolliding random walks; lozenge tilings}

\Classification{06C05; 05E05; 05A30}

\section{Introduction}\label{sec:intro}

\subsection{Overview}

The Dyson Brownian motion \cite{dyson1962brownian}
is a continuous stochastic dynamics of $N$
particles on the one-dimensional line $\mathbb{R}$.
The particles evolve according to independent
Brownian motions which are conditioned to never collide.
The noncolliding property may be also modeled as
Coulomb repelling.
The Dyson Brownian motion
arises
as the dynamics of eigenvalues from the standard Brownian motion
on the space of complex Hermitian matrices. As such,
it has been heavily utilized towards universality
results for random matrix spectra
\cite{AndersonGuionnetZeitouniBook,LY_RMT_Bull2011,Johansson2001Universality}.

Within integrable probability, a number of discrete deformations of the
Dyson Brownian motion were introduced, starting from
noncolliding Poisson and Bernoulli random walks
\cite{konig2002non}
(based on a classical formula of \cite{KMG59-Coincidence})
and
followed by
their Macdonald deformation depending on two parameters $(q,t)$ which is
defined in~\cite{BorodinCorwin2011Macdonald}.
A notable special case of the latter
considered in~\cite{GorinShkolnikov2014}
is the Jack limit $t=q^{\beta/2}\to 1$, where $\beta>0$ is the
beta parameter from random matrix theory. In the Jack limit,
one obtains $\beta$-noncolliding Poisson random walks
(further studied in~\cite{huang2021beta}),
and also their multilevel versions.
A scaling limit of the latter leads to the multilevel
Dyson Brownian motion with the general $\beta$ parameter.

Each of the known discrete deformations of the Dyson Brownian motion
is powered by the Cauchy summation
identity for some family of symmetric polynomials
$\left\{ P_\lambda \right\}$
such as
Schur (for noncolliding Poisson and Bernoulli walks),
Jack, or Macdonald polynomials.
Here $\left\{ P_\lambda \right\}$ is one of these families
of polynomials
in $N$ variables. The $P_\lambda$'s
form a linear basis in the space
of symmetric polynomials in $N$ variables as $\lambda$ runs
over
partitions $\lambda=(\lambda_1\ge \dots\ge \lambda_N\ge0 )$,
$\lambda_i\in \mathbb{Z}$, with~$N$ parts.

The Cauchy identity is a fundamental property of
many families of symmetric polynomials, and is
closely tied to their orthogonality with respect
to a suitable inner product.
It provides
a~product-form expression
for the sum $\sum_\lambda b_\lambda
P_\lambda(x_1,\dots,x_N )
P_\lambda(y_1,\dots,y_N )$,
where $b_\lambda$ are certain explicit coefficients.
In stochastic dynamics of $N$ noncolliding particles,
this Cauchy identity implies the normalization to one
property of the transition probability, which also involves an
$N$-fold summation over partitions $\lambda$.

Along with the Cauchy identity,
most families of symmetric polynomials
satisfy a branching rule.
This identity
expresses $P_\lambda(x_1,\dots,x_N )$
in $N$ variables as a nonnegative linear combination
of polynomials $P_\mu(x_1,\dots,x_{N-1} )$
with $N-1$ variables, where the sum runs over $\mu$.
For particular symmetric polynomials, such an expansion has clear
representation-theoretic meaning.
For example,
for Schur polynomials the branching rule
is behind the decomposition of a given irreducible representation
of the unitary group $U(N)$
when restricted to the subgroup $U(N-1)$.

Note that the branching rule is often dual to the Pieri rule
expressing the product $f\hspace{1pt} P_\lambda$ (for a special choice of $f$
like $x_1+\dots+x_N$)
as a linear combination of $P_\nu$'s in the same number of variables.
In the present paper we do not explicitly use this duality,
and keep the branching rule perspective.

The goal of the present work is to
construct and
explore noncolliding random walks arising
from the branching rule instead of the Cauchy identity.
We start at the level of Macdonald polynomials
with parameters $(q,t)\in(0,1)^2$, and
take a limit
of the branching rule under the principal specialization
$(x_1,\dots,x_N )=\big(1,t,\dots,t^{N-1} \big)$
as the number of variables $N$ goes to infinity.
\mbox{Using} the resulting summation identity
(formulated in Theorem~\ref{thm:Macdonald_noncolliding_intro}
later in the introduction),
we define a
new
discrete-time Markov process $\Upsilon_m$ of $m$ distinct ordered
particles in $\mathbb{Z}_{\ge0}$
with negative drift and absorbing wall
at zero (where $m\in \mathbb{Z}_{\ge1}$
is assumed fixed).
The presence of the wall means that the process
almost surely
reaches its only absorbing state $(m-1,m-2,\dots,2,1,0 )$.

Trajectories of $\Upsilon_m$ may be identified with
lozenge tilings or plane partitions with certain
explicit boundary conditions depending on
the initial configuration in
$\Upsilon_m$.
We show that the probability measure on plane partitions coming from
$\Upsilon_m$ has a Gibbs characterization via the so-called
Boltzmann factors
which are ratios of
probability weights of two configurations
differing by an elementary transformation.
We explicitly compute these Boltzmann factors
in the general Macdonald case.
In the particular case $t=q$,
the Gibbs probability weight
of a plane partition is proportional simply to
$q^{\mathsf{vol}}$, where $\mathsf{vol}$ is the sum of the
entries of the plane partition.
See Section~\ref{sec:plane_partitions} in the text for details.

We consider a number of degenerations of our process
$\Upsilon_m$ leading to known deformations
of the Dyson Brownian motion mentioned above.
All these degenerations correspond to specializing the
parameters $(q,t)$ in such a way that the Macdonald
polynomials turn into another well-known family of symmetric polynomials:
\begin{itemize}\itemsep=0pt
	\item (Schur polynomials)
		Setting $t=q$, we get a simpler
		Markov process
		$\Upsilon_m^{\mathrm{Schur}}$
		of $m$ noncolliding particles
		on $\mathbb{Z}_{\ge0}$
		with an absorbing wall at $0$.
		To the best of the author's knowledge,
		this process and the underlying normalization identity
		(stating that the quantities in \eqref{eq:Upsilon_Schur} below
		sum to~$1$) are also new.
		The process $\Upsilon_m^{\mathrm{Schur}}$
		looks similar to the translation invariant
		$q$-noncolliding random walks on $\mathbb{Z}$
		introduced
		and studied in
		\cite{BG2011non}.
		The
		normalization of transition probabilities
		in the latter process
		can be traced back to the Cauchy
		identity. However, it does not seem that our process
		$\Upsilon_m^{\mathrm{Schur}}$ can be scaled to that of
		\cite{BG2011non}.

	\item (Jack polynomials)
		Take $t=q^{\beta/2}\to 1$ (where $\beta>0$ is the parameter
		coming from random matrix theory),
		and simultaneously
		scale the coordinates of the process $\Upsilon_m$
		away from $0$. In this way we get a dynamics of $m$
		particles on $\mathbb{Z}$ which is invariant under
		space translations of the particles.
		This dynamics is closely related to the
		$\beta$-noncolliding Poisson random walks
		studied in \cite{GorinShkolnikov2014,huang2021beta}.
		Thus, we see that our new Macdonald noncolliding walks~$\Upsilon_m$ generalize all known noncolliding processes
		at the Jack level (with general random matrix~$\beta$ parameter).
		In the particular case $\beta=2$,
		we recover the Bernoulli and Poisson
		walks conditioned to never collide
		which were studied in~\cite{konig2002non}. Under Brownian scaling,
		it is known that the latter random walks turn into the
		classical Dyson Brownian motion coming from Hermitian
		random matrices.

	\item (Hall--Littlewood polynomials)
		Setting
		$q=0$, further sending $t\to1$ and taking a Poisson-type limit from
		discrete time to continuous,
		we arrive at a new particle system
		$\Upsilon_m^{\mathrm{cont}}$
		on $\mathbb{Z}_{\ge0}$
		with an absorbing wall at zero which evolves as follows.
		To each particle $\mathsf{x}_1>\dots>\mathsf{x}_m\ge0 $
		we assign an independent exponential clock
		of rate $i(\mathsf{x}_i-\mathsf{x}_{i+1}-1)$,
		where, by agreement, $\mathsf{x}_{m+1}=-1$.
		When the clock of $\mathsf{x}_{i}$ rings,
		we additionally
		select an index $j \in \{ 1,\dots,i \}$
		uniformly at random,
		and all the particles
		$\mathsf{x}_i,\mathsf{x}_{i-1},\dots,\mathsf{x}_j $
		simultaneously
		jump to the left by $1$.
		The process $\Upsilon_m^{\mathrm{cont}}$
		almost surely reaches its absorbing state
		$(m-1,\dots,1,0)$.
		A more detailed investigation of this particle system will be
		performed elsewhere.
\end{itemize}

In the next
Section~\ref{sub:intro_details}
we describe in detail our most general Markov processes
$\Upsilon_m$ arising at the Macdonald level.

\subsection{Macdonald noncolliding walks}\label{sub:intro_details}

Throughout the paper we assume that
$(q,t)$ are real numbers belonging to $(0,1)$.
We need some notation.
Recall that the $q$-Pochhammer symbols are given by
\begin{equation*}
	(z;q)_k\coloneqq (1-z)(1-zq)\cdots\big(1-zq^{k-1}\big),
	\qquad
	k\in \mathbb{Z}_{\ge0},
\end{equation*}
and $(z;q)_\infty\coloneqq\prod_{i=0}^{\infty}(1-zq^i)$ is a convergent infinite
product because $|q|<1$.

For $\vec{\mathsf{x}}=(\mathsf{x}_1>\dots>\mathsf{x}_m\ge0 )$,
denote the $(t,q)$-deformed Vandermonde product by
\begin{equation}
	\label{eq:V_tq}
	V_{t,q}(\vec{\mathsf{x}}):=
	\prod_{1\le i<j\le m}
	\frac{\left( q^{j-i-1}t^{\mathsf{x}_i-\mathsf{x}_j-j+i+1};t \right)_{\infty}}
	{(q^{j-i} t^{\mathsf{x}_i-\mathsf{x}_j-j+i+1};t)_{\infty}}.
\end{equation}
When $t=q\to 1$, $V_{t,q}$ turns (after rescaling by a suitable power of $\log q$)
into the
usual Vandermonde
$V(\vec{\mathsf{x}})=\prod_{1\le i<j \le m}(\mathsf{x}_i-
\mathsf{x}_j)$. Moreover,
when $\mathsf{x}_i=\mathsf{x}_{i+1}$ for some $i$,
one readily sees that $V_{t,q}(\vec{\mathsf{x}})$
vanishes.

Let $\mathsf{y}=(\mathsf{y}_1>\dots>\mathsf{y}_m\ge0 )$
be such that
\begin{equation}\label{eq:y_differs_from_x}
	\mathsf{y}_i-\mathsf{x}_i\in \{ -1,0 \}
	\qquad
	\text{for all $1\le i\le m$}.
\end{equation}
Now we can define the main object of the present paper:
\begin{gather}
\Upsilon_m(\vec{\mathsf{x}},\vec{\mathsf{y}})
		\coloneqq
		t^{-\binom m2}	
		\hspace{1pt}
		\frac
		{V_{t,q}(\vec{\mathsf{y}})}
		{V_{t,q}(\vec{\mathsf{x}})}
		\prod_{\substack{1\le i<j\le m
		\\
		\mathsf{y}_{i}=\mathsf{x}_{i}
		,\,
		\mathsf{y}_{j}=\mathsf{x}_{j}-1}}
		\frac{\big(1-t^{i-j+\mathsf{x}_{i}-\mathsf{x}_{j}+1}
		q^{j-i-1}\big)
		\big(1-t^{i-j+\mathsf{x}_{i}-\mathsf{x}_{j}} q^{j-i+1}\big)}
		{\big(1-t^{i-j+\mathsf{x}_{i}-\mathsf{x}_{j}+1} q^{j-i}\big)
		\big(1-t^{i-j+\mathsf{x}_{i}-\mathsf{x}_{j}} q^{j-i}\big)}
		\nonumber\\
\hphantom{\Upsilon_m(\vec{\mathsf{x}},\vec{\mathsf{y}}) \coloneqq}{}
\times\prod_{i\colon \mathsf{y}_i=\mathsf{x}_i}
		t^{\mathsf{x}_i}
		\prod_{i\colon \mathsf{y}_i=\mathsf{x}_i-1}
		\big(t^{m-i}-q^{m-i}t^{\mathsf{x}_i}\big)		.\label{eq:intro_Upsilon}
\end{gather}
From \eqref{eq:y_differs_from_x}
one readily sees that the infinite products in
${V_{t,q}(\vec{\mathsf{y}})}/{V_{t,q}(\vec{\mathsf{x}})}$
cancel out in such a way that
\eqref{eq:intro_Upsilon} is always a rational function of $q$, $t$.
Moreover, for $0<q,t<1$ the quantities \eqref{eq:intro_Upsilon}
are nonnegative.
One of our main results is the sum-to-one identity for the
$\Upsilon_m$'s:
\begin{Theorem}
	\label{thm:Macdonald_noncolliding_intro}
	With the above notation,
	for any $\vec{\mathsf{x}}=(\mathsf{x}_1>\dots>\mathsf{x}_m\ge0 )$
	we have
	\begin{equation*}
		\sum_{
			\substack{\vec{\mathsf{y}}=
			(\mathsf{y}_1>\dots>\mathsf{y}_m\ge0 )
			\\
			\mathsf{y}_i-\mathsf{x}_i\in \left\{ -1,0 \right\}
			\
			\text{for all $1\le i\le m$}
			}
		}
		\Upsilon_m(\vec{\mathsf{x}},\vec{\mathsf{y}})=1.
	\end{equation*}
\end{Theorem}
Theorem~\ref{thm:Macdonald_noncolliding_intro} implies that
the quantities
$\Upsilon_m(\vec{\mathsf{x}},\vec{\mathsf{y}})$
may be viewed as
transition probabilities of a discrete time Markov chain
of $m$ ordered distinct particles
on $\mathbb{Z}_{\ge0}$ in which at
each step, each particle either stays, or moves to the left by $1$.
Eventually with probability $1$ the chain reaches the
absorbing state $(m-1,m-2,\dots,2,1,0 )$,
see Proposition~\ref{prop:absorbing_wall}
below.
We call the Markov chain $\Upsilon_m$ the \emph{Macdonald
noncolliding walks} with an absorbing wall at zero.

We prove Theorem~\ref{thm:Macdonald_noncolliding_intro}
in Section~\ref{sec:limit_transition}
by obtaining the transition probabilities
$\Upsilon_m(\vec{\mathsf{x}},\vec{\mathsf{y}})$
as a limit of certain
ratios of Macdonald polynomials
evaluated at the principal specializations
$\big(1,t,t^2,\dots,t^{N-1} \big)$, as the number of variables
goes to infinity. The fact that before the limit
these ratios sum to one is equivalent to the branching
rule for the Macdonald polynomials.

\subsection{Outline}

In Section~\ref{sec:macdonald_polynomials} we review the definition
of Macdonald symmetric polynomials together
with all the required formulas.
In Section~\ref{sec:limit_transition} we perform the main limit
transition, and obtain the Macdonald noncolliding walks
$\Upsilon_m$.
In Section~\ref{sec:properties_of_our_walks}
we consider various degenerations of our dynamics
when the Macdonald parameters $(q,t)$ are specialized in a certain way.
More precisely, we look at the dynamics at $t=q$ (when the Macdonald polynomials
reduce to the Schur polynomials), as $q=t^\alpha\to1$
(reduction to the Jack polynomials),
and as $q=0$ (when Macdonald polynomials become the Hall--Littlewood polynomials).
Moreover, in the latter case we see that sending $t\to1$ leads
to a new continuous time Markov chain on $m$ particles with
inhomogeneous jump rates.
In Section~\ref{sec:plane_partitions} we
give a Gibbs characterization of the probability measure
on the space of trajectories of our noncolliding walks
by means of the so-called Boltzmann factors which are ratios of
probability weights of two trajectories differing by an elementary transformation.

\section{Review of Macdonald polynomials}\label{sec:macdonald_polynomials}

Here we collect the necessary notation and results
around Macdonald symmetric polynomials.
We follow \cite[Chapter VI]{Macdonald1995}.

\subsection{Definition}\label{sub:Macdonald_polynomials_definition}

Let $N\ge1$.
Macdonald symmetric polynomials $P_\lambda$ in $N$ variables $x_1,\dots,x_N$
are indexed by partitions
$\lambda=(\lambda_1\ge \dots\ge \lambda_N \ge0)$, $\lambda_i\in \mathbb{Z}$,
with $N$ parts. Denote the set of these partitions by~$\mathbb{Y}(N)$.
The $P_\lambda$'s
depend on two parameters $q,t\in[0,1)$. For each fixed
$(q,t)$, they form a~basis in the space of symmetric
polynomials in~$N$ variables when $\lambda$ runs over~$\mathbb{Y}(N)$.
The shortest definition of the
$P_\lambda$'s is through the first Macdonald $q$-difference operator
acting in the~$x_i$'s:
\begin{gather*}
	D_1\coloneqq \sum_{i=1}^{N}
	\biggl( \prod_{1\le j\le N\colon j\ne i}\frac{x_j-tx_i}{x_j-x_i} \biggr)
	T_{q;i},\\ T_{q;i}f(x_1,\dots,x_N )\coloneqq f(x_1,\dots,x_{i-1},qx_i,x_{i+1},\dots,x_N ).
\end{gather*}
The operator $D_1$ preserves the space of symmetric
polynomials in $x_1,\dots,x_N $, and its eigenfunctions
are the Macdonald polynomials
\begin{gather}
	D_1P_\lambda(x_1,\dots,x_N\mid q,t )=
	\bigl( q^{\lambda_1}t^{N-1}+q^{\lambda_2}t^{N-2}+\dots+q^{\lambda_N}
	\bigr)
	P_\lambda(x_1,\dots,x_N\mid q,t ),\label{eq:Macdonald_eigenrelation}
\end{gather}
with $\lambda\in \mathbb{Y}(N)$.
For generic $(q,t)$, the eigenvalues in~\eqref{eq:Macdonald_eigenrelation}
for different $\lambda$ are distinct, so
the $P_\lambda$'s are determined uniquely up to normalization.
The normalization is specified by
\begin{equation*}
	P_\lambda(x_1,\dots,x_N\mid q,t )=
	x_1^{\lambda_1}x_2^{\lambda_2}\cdots x_N^{\lambda_N} +
	\text{lexicographically lower terms},
\end{equation*}
where the lower terms depend on $q$, $t$.

In the case $q=t$, the polynomials $P_\lambda$ reduce to the
well-known Schur symmetric polynomials~$s_\lambda$, which admit
the following explicit determinantal formula (which does not
depend on the choice of $q$):
\begin{equation}
	\label{eq:Schur}
	P_\lambda(x_1,\dots,x_N\mid q,q )=s_\lambda(x_1,\dots,x_N )=
	\frac{\det\bigl[ x_i^{\lambda_j+N-j} \bigr]_{i,j=1}^{N}}{\det
	\bigl[ x_i^{N-j} \bigr]_{i,j=1}^N}.
\end{equation}
For generic $(q,t)$, there are no known formulas for Macdonald polynomials which are as compact as~\eqref{eq:Schur}.

\subsection{Principal specialization}\label{sub:principal_specialization}

\begin{figure}[htpb]
	\centering
	\includegraphics[width=.28\textwidth]{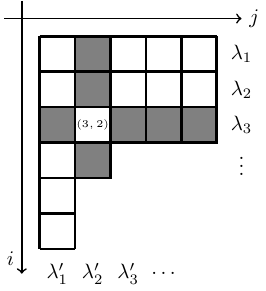}
	\caption{Young diagram $\lambda=(5,5,5,2,1,1)$. Highlighted
		are the arm, leg, coarm, and coleg of the box $\square=(3,2)$.
		We have $|\lambda|=19$ and
		$a(\square)=3$, $l(\square)=1$, $a'(\square)=1$, $l'(\square)=2$.}
	\label{fig:Young_d}
\end{figure}

When the variables $x_i$ are specialized to a finite geometric
progression with the ratio $t$ (the second Macdonald parameter),
the polynomials $P_\lambda$ admit an explicit product formula in terms
of the Young diagram corresponding to the partition $\lambda$.
Recall that the Young diagram $\lambda$ is a~collection of $1\times 1$ boxes in the plane
with $\lambda_i$ boxes in row $i$, see Figure~\ref{fig:Young_d}
for an illustration.
The principal specialization takes the form
\cite[formulas~(VI.6.11) and (VI.6.11$'$)]{Macdonald1995}:
\begin{equation}	\label{eq:principal_spec}
	P_\lambda\big(1,t,\dots,t^{N-1} \mid q,t \big)=
	t^{n(\lambda)}\prod_{\square\in \lambda}
	\frac{1-q^{a'(\square)}t^{N-l'(\square)}}
	{1-q^{a(\square)}t^{l(\square)+1}},
\end{equation}
where the product is over all boxes of the Young diagram $\lambda$,
\begin{equation}	\label{eq:n_lambda_notation}
	n(\lambda):=\sum_i(i-1)\lambda_i,
\end{equation}
and the arms, legs, coarms, colegs of the box $\square=(i,j)$
are defined, respectively, as
\begin{equation}	\label{eq:arms_legs}
	a(\square)=\lambda_i-j,\qquad l(\square)=\lambda_j'-i,\qquad
	a'(\square)=j-1,\qquad l'(\square)=i-1.
\end{equation}
Here $\lambda_j'$ are the column lengths in the Young
diagram, see Figure~\ref{fig:Young_d}.

\subsection{Branching}\label{sub:branching}

Let us first recall the Pieri coefficients for Macdonald polynomials
\cite[formula~(VI.6.24.ii)]{Macdonald1995},
\cite[formula~(2.11)]{BorodinCorwin2011Macdonald}.
They depend on a pair of partitions $\mu$, $\lambda$
which \emph{interlace}, namely,
\begin{equation*}%\label{eq:interlace}
	\lambda_1\ge \mu_1\ge \lambda_2\ge \mu_2\ge \cdots
\end{equation*}
(notation $\mu\prec \lambda$), and are defined as
\begin{equation}
	\label{eq:psi_Pieri}
	\psi_{\lambda/\mu}
	=
	\psi_{\lambda/\mu}(q,t)
	\coloneqq
	\prod_{1\le i<j\le \ell(\mu)}
	\frac{\mathsf{f}\big(q^{\mu_i-\mu_j}t^{j-i}\big)
	\hspace{1pt}\mathsf{f}\big(q^{\lambda_i-\lambda_{j+1}}t^{j-i}\big)}
	{\mathsf{f}\big(q^{\lambda_i-\mu_j}t^{j-i}\big)
	\hspace{1pt}\mathsf{f}\big(q^{\mu_i-\lambda_{j+1}}t^{j-i}\big)},
	\qquad
	\mathsf{f}(u):=\frac{(tu;q)_{\infty}}{(qu;q)_{\infty}}.
\end{equation}
Here and below by $\ell(\mu)$ we denote the number of
parts in $\mu$ which are strictly positive,
that is, $\mu_{\ell(\mu)}>0$, $\mu_{\ell(\mu)+1}=0$ (when $\mu=(0,0,\dots )$, we set $\ell(\mu)=0$).
Let also $|\lambda|=\lambda_1+\dots+\lambda_{\ell(\lambda)}$
denote the number of boxes in the Young diagram~$\lambda$.

\begin{Proposition}[{\cite[formula~(VI.7.13$'$)]{Macdonald1995}}]
	\label{prop:branching}
	Let $\lambda\in \mathbb{Y}(N)$.
	We have
	\begin{equation}
		\label{eq:branching}
		P_\lambda(x_1,\dots,x_N \mid q,t)
		=
		\sum_{\mu\colon\mu\prec\lambda}
		\psi_{\lambda/\mu}(q,t)\hspace{1pt} P_\mu(x_1,\dots,x_{N-1} \mid q,t )\hspace{1pt}
		x_N^{|\lambda|-|\mu|},
	\end{equation}
	where the sum is over
	$\mu\in \mathbb{Y}(N-1)$ which interlace with $\lambda$.
\end{Proposition}

Note that for $q,t\in [0,1)$ the coefficients $\psi_{\lambda/\mu}$
\eqref{eq:psi_Pieri} are all nonnegative.
Together with Proposition~\ref{prop:branching} this implies:
\begin{Corollary}
	\label{cor:nonnegative_spec}
	Let $q,t\in[0,1)$. Specializing the variables
	$x_1,\dots,x_N $ into
	nonnegative real numbers
	makes the Macdonald polynomial
	$P_\lambda(x_1,\dots,x_N\mid q,t)$ nonnegative.
\end{Corollary}

More generally, for $N\ge K\ge1$
define the skew Macdonald polynomials
$P_{\lambda/\mu}$
as the coefficients
in the expansion:
\begin{equation*}
	P_\lambda(x_1,\dots,x_N \mid q,t)=
	\sum_{\mu\in \mathbb{Y}(K)}
	P_\mu(x_1,\dots,x_K \mid q,t )
	P_{\lambda/\mu}(x_{K+1},\dots,x_{N-1},x_N \mid q,t ).
\end{equation*}
Here we use the fact that
the polynomials
$P_\mu(x_1,\dots,x_K )$
form a basis in the space of symmetric functions in $K$ variables, and
expand
$P_\lambda(x_1,\dots,x_N )$
in this basis.
The Pieri coefficient is related to the skew Macdonald polynomial
in one variable as follows:
\begin{equation*}%\label{eq:P_skew_single_variable}
	P_{\lambda/\mu}(x_1\mid q,t)=
	\begin{cases}
		\psi_{\lambda/\mu}(q,t)\hspace{1pt} x_1^{|\lambda|-|\mu|},&\mu\prec \lambda,\\
		0,&\text{otherwise}.
	\end{cases}
\end{equation*}

Later we will also use the dual Pieri coefficients
$\psi_{\lambda/\mu}'$ which are defined as
(see \cite[formula~(VI.6.24.iv)]{Macdonald1995} and
\cite[formula~(2.12)]{BorodinCorwin2011Macdonald})
\begin{equation}
	\label{eq:psi_prime_Pieri}
	\psi_{\lambda/\mu}'=
	\psi_{\lambda/\mu}'(q,t)
	\coloneqq
	\prod_{\substack{i<j \\ \lambda_{i}=\mu_{i},\, \lambda_{j}=\mu_{j}+1}}
	\frac{\big(1-q^{\mu_{i}-\mu_{j}} t^{j-i-1}\big)\big(1-q^{\lambda_{i}-\lambda_{j}} t^{j-i+1}\big)}
	{\big(1-q^{\mu_{i}-\mu_{j}} t^{j-i}\big)\big(1-q^{\lambda_{i}-\lambda_{j}} t^{j-i}\big)}.
\end{equation}
for partitions $\lambda$, $\mu$ whose column lengths interlace,
that is,
$\lambda_1'\ge \mu_1'\ge \lambda_2'\ge \mu_2'\ge \cdots $.
We have
\begin{equation}
	\label{eq:P_skew_single_variable_dual}
	P_{\lambda'/\mu'}(x_1 \mid t,q)=
	\begin{cases}
		\psi'_{\lambda/\mu}(q,t)\hspace{1pt} x_1^{|\lambda|-|\mu|}
		,&\mu'\prec \lambda',\\
		0
		,&\text{otherwise}.
	\end{cases}
\end{equation}

\subsection{Cauchy identity}\label{sub:Cauchy_and_Dyson1}

Along with the branching rule, another fundamental
identity for Macdonald polynomials is the Cauchy identity.
Here we present its dual version, see \cite[Chapter~VI.4]{Macdonald1995} for the
usual version.

\begin{Proposition}[{\cite[formula~(VI.5.4) and Chapter~VI.7]{Macdonald1995}}]
	\label{prop:dual_Cauchy_identity}
	Let $N,M\ge 1$ be fixed. We have
	\begin{equation*}%\label{eq:dual_Cauchy}
		\sum_{\lambda\in \mathbb{Y}(N)\colon \lambda_1\le M}
		P_\lambda(x_1,\dots,x_N \mid q,t)
		\hspace{1pt}
		P_{\lambda'}(y_1,\dots,y_M\mid t,q )=
		\prod_{i=1}^{N}\prod_{j=1}^{M}(1+x_iy_j).
	\end{equation*}
	Moreover, for any $\mu\in \mathbb{Y}(N)$ we have
	the following particular case of
	the skew Cauchy identity:
	\begin{equation}		\label{eq:dual_Cauchy_skew}
		P_\mu(x_1,\dots,x_N\mid q,t )
		\prod_{i=1}^{N}(1+x_iy)=
		\sum_{\lambda\in \mathbb{Y}(N)\colon \mu'\prec \lambda'}
		P_{\lambda'/\mu'}(y\mid t,q)
		\hspace{1pt}
		P_\lambda(x_1,\dots,x_N\mid q,t ).
	\end{equation}
\end{Proposition}
Identities~\eqref{eq:branching} and~\eqref{eq:dual_Cauchy_skew} look very similar, but note that
in the former we sum over the smaller partition, while in the latter one we sum over the larger partition.

\subsection{Markov kernels from Macdonald polynomials}\label{eq:Markov_kernels_from_Cauchy}

Let us first recall a general definition from
\cite[Chapter~7]{borodin2016representations}.
Let $\mathfrak{X}$, $\mathfrak{Y}$ be
finite or countable sets. By a Markov kernel
(or a link) from $\mathfrak{X}$ to $\mathfrak{Y}$,
we mean a function $P$ on $\mathfrak{X}\times \mathfrak{Y}$
such that
$P(x,y)\in[0,1]$ for all $x\in \mathfrak{X}$, $y\in \mathfrak{Y}$,
and
\begin{equation*}
	\sum_{y\in \mathfrak{Y}}P(x,y)=1
	\qquad
	\text{for all $x\in \mathfrak{X}$}.
\end{equation*}
We adopt the notation $P\colon \mathfrak{X}\dashrightarrow \mathfrak{Y}$.

Normalizing identities
\eqref{eq:branching} and \eqref{eq:dual_Cauchy_skew}
with nonnegative variables $x_i$ and $y$
(following
\cite{Borodin2010Schur,BorFerr2008DF} in the Schur case and~\cite{BorodinCorwin2011Macdonald}
in the general Macdonald case)
leads to the following two families of links:
\begin{align}
\Lambda^{N}_{N-1}\colon \ & \mathbb{Y}(N)\dashrightarrow \mathbb{Y}(N-1),
		\nonumber\\
		& \Lambda^{N}_{N-1}(\lambda,\mu)\coloneqq
		\psi_{\lambda/\mu}(q,t)\hspace{1pt} x_N^{|\lambda|-|\mu|}\hspace{1pt}
		\frac{
		P_\mu(x_1,\dots,x_{N-1} \mid q,t )
		}
		{P_\lambda(x_1,\dots,x_N \mid q,t)}
		\,\mathbf{1}_{\mu\prec\lambda},\label{eq:traditional_links_Lambda}\\
Q_N\colon \ & \mathbb{Y}(N)\dashrightarrow \mathbb{Y}(N),
		\nonumber\\
		& Q_N(\lambda,\nu)\coloneqq
		\frac{\psi'_{\nu/\lambda}(q,t)
		\hspace{1pt}
		y^{|\nu|-|\lambda|}
		}
		{\prod_{i=1}^{N}(1+x_i y)}
		\frac
		{
		P_\nu(x_1,\dots,x_N\mid q,t )}
		{P_\lambda(x_1,\dots,x_N\mid q,t )}
		\,\mathbf{1}_{\lambda'\prec\nu'}		.\label{eq:traditional_links_Q}
\end{align}
Here and
throughout the paper by $\mathbf{1}_{A}$ we denote
the indicator of the event (or condition) $A$.
For the link \eqref{eq:traditional_links_Q}
we also used \eqref{eq:P_skew_single_variable_dual}.

The links \eqref{eq:traditional_links_Lambda}--\eqref{eq:traditional_links_Q}
satisfy the following intertwining
relation:
\begin{gather}
	Q_N
	\Lambda^{N}_{N-1}
	=
	\Lambda^{N}_{N-1}
	Q_{N-1}
	,\nonumber\\
	\sum_{\nu\in \mathbb{Y}(N)}
	Q_N(\lambda,\nu)
	\Lambda^{N}_{N-1}(\nu,\varkappa)
	=
	\sum_{\mu\in \mathbb{Y}(N-1)}
	\Lambda^{N}_{N-1}(\lambda,\mu)
	Q_{N-1}(\mu,\varkappa)
	,\label{eq:traditional_intertwining}
\end{gather}
where the second identity holds for all $\lambda\in \mathbb{Y}(N)$,
$\varkappa\in \mathbb{Y}(N-1)$,
and is simply a more detailed rewriting of the first one.
Intertwining relation~\eqref{eq:traditional_intertwining}
follows from the skew Cauchy identity,
for example, see \cite[Proposition~2.3.1]{BorodinCorwin2011Macdonald}.

The Markov chain on $\mathbb{Y}(N)$ defined by the operator
$Q_N$ is traditionally viewed
as the
Macdonald deformation of the
Dyson Brownian motion, see, for example,
\cite{BorFerr2008DF, GorinShkolnikov2012}
for the Schur $q=t$ case, and also
\cite{GorinShkolnikov2014} for the general $\beta$ version
based on Jack symmetric polynomials.
For the classical $\beta=2$ Dyson Brownian motion (coming from the
Gaussian unitary ensemble of random matrices),
intertwining relations were investigated in \cite{warren2005dyson}.

In this paper we take a limit of the links
$\Lambda^{N}_{N-1}$ as $N\to\infty$ to construct
a new Markov process~$\Upsilon_m$ of~$m$ noncolliding particles
depending on the Macdonald parameters~$(q,t)$.
This process, too,
may be viewed as another Macdonald deformation
of the Dyson Brownian motion (in particular,
our process admits a diffusive scaling
to the Dyson Brownian motion).
We also note that in the limit we consider, the matrix elements of the
operators $Q_N$ corresponding to our scaling tend to zero. Thus, it is not clear
whether the limiting Markov chains coming from~$\Lambda^{N}_{N-1}$ admit any intertwining relation like \eqref{eq:traditional_intertwining}.
In contrast, they are going to be consistent for
different numbers of particles, see
Proposition~\ref{prop:consistency} below.

\section{Limit transition to noncolliding walks}\label{sec:limit_transition}

In this section we perform the limit transition as $N\to+\infty$ in the
Markov kernels $\Lambda^{N}_{N-1}$ \eqref{eq:traditional_links_Lambda}
under the principal specialization
\eqref{eq:principal_spec}, and prove
Theorem~\ref{thm:Macdonald_noncolliding_intro}.

\subsection{Setup}\label{sub:limit_transition_setup}

\begin{figure}[htpb]
	\centering
	\includegraphics[width=.8\textwidth]{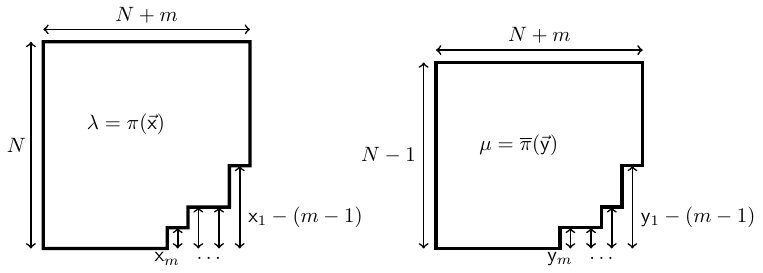}
	\caption{Left: Young diagram of $\pi(\vec{\mathsf{x}})$ defined by
	\eqref{eq:pi_maps_1}--\eqref{eq:pi_maps_2}.
	The column lengths of the defect are $\mathsf{x}_j-(m-j)$.
	Right: Young diagram of $\overline\pi(\vec{\mathsf{y}})$.
	Observe that
	$\lambda=\pi(\vec{\mathsf{x}})$ and
	$\overline\pi(\vec{\mathsf{y}})$ differ only by
	adding the first row.}	\label{fig:pi_of_x}
\end{figure}

Denote by $\mathbb{W}_m$ the space of $m$-particle configurations
in $\mathbb{Z}_{\ge0}$, that is,
\begin{equation}
	\label{eq:W_m_space}
	\mathbb{W}_m\coloneqq
	\{
		\vec{\mathsf{x}}=(\mathsf{x}_1>\mathsf{x}_2>\dots>\mathsf{x}_m\ge0 )
	\}\subset \mathbb{Z}_{\ge0}^{m}.
\end{equation}
By $\mathbb{W}_m(N)$ denote the finite subset of $\mathbb{W}_m$ determined by
the condition
$\mathsf{x}_1\le N+m-2$. Define the injective maps
\begin{equation}
	\label{eq:pi_maps_1}
	\pi\colon \ \mathbb{W}_m(N)\to \mathbb{Y}(N), \qquad
	\overline\pi\colon \ \mathbb{W}_m(N)\to \mathbb{Y}(N-1),
\end{equation}
as follows. If $\lambda=\pi(\vec{\mathsf{x}})\in \mathbb{Y}(N)$
and $\mu=\overline\pi(\vec {\mathsf{y}})\in \mathbb{Y}(N-1)$, then
\begin{gather}
 \{ \lambda_1-1,\lambda_2-2,\dots,\lambda_N-N \}=
 \{ 0,1,2,\dots,N+m-1 \}\setminus
		\{\mathsf{x}_1,\dots, \mathsf{x}_m \},\nonumber\\
		 \{ \mu_1-1,\mu_2-2,\dots,\mu_{N-1}-(N-1) \}=
		 \{ 1,2,\dots,N+m-1 \}\setminus
		\{\mathsf{y}_1+1,\dots, \mathsf{y}_m+1 \}.\!\!\!\label{eq:pi_maps_2}
	\end{gather}
For fixed $\vec{\mathsf{x}}$ and growing $N$,
almost all parts of $\lambda=\pi(\vec{\mathsf{x}})$
are equal to $N+m$, and there is a defect in a few last parts of
$\lambda$. The columns of this defect are
encoded through
$\vec{\mathsf{x}}$. A similar description holds for
$\mu=\overline\pi(\vec{\mathsf{y}})$.
In multiplicative notation for partitions, we have
\begin{gather}
		\lambda=\pi(\vec{\mathsf{x}})=
		(N+m)^{N-\mathsf{x}_1+m-1}(N+m-1)^{\mathsf{x}_1-\mathsf{x}_2-1}
		\cdots (N+1)^{\mathsf{x}_{m-1}-\mathsf{x}_m-1} N^{\mathsf{x}_m},
\nonumber\\
		\mu=\overline\pi(\vec{\mathsf{y}})=
		(N+m)^{N-\mathsf{y}_1+m-2}(N+m-1)^{\mathsf{y}_1-\mathsf{y}_2-1}
		\cdots (N+1)^{\mathsf{y}_{m-1}-\mathsf{y}_m-1}
		N^{\mathsf{y}_m}.\label{eq:pi_maps_multiplicative_notation}
	\end{gather}
See
Figure~\ref{fig:pi_of_x} for an illustration.

\begin{Lemma}	\label{lemma:interlacing}
	For any $\vec{\mathsf{x}},\vec{\mathsf{y}}\in \mathbb{W}_m(N)$,
	we have $\overline\pi(\vec{\mathsf{y}})\prec \pi(\vec{\mathsf{x}})$
	if and only if $\mathsf{y}_i=\mathsf{x}_i$ or $\mathsf{y}_i=\mathsf{x}_i-1$
	for all $i=1,\dots,m $.
\end{Lemma}
\begin{proof}
	Straightforward verification.
\end{proof}

In the rest of this section we
fix arbitrary $\vec{\mathsf{x}},\vec{\mathsf{y}}\in \mathbb{W}_m$ and
compute the limit of
$\Lambda^{N}_{N-1}(\pi(\vec{\mathsf{x}}),\overline\pi(\vec{\mathsf{y}}))$
\eqref{eq:traditional_links_Lambda}
under principal
specialization $x_j=t^{j-1}$
as $N\to+\infty$. Clearly, for sufficiently large $N$ the
partitions
$\pi(\vec{\mathsf{x}}),\overline\pi(\vec{\mathsf{y}})$
are well-defined.
We will show that this limit is equal
to the Markov kernel $\Upsilon_m(\vec{\mathsf{x}},\vec{\mathsf{y}})$
defined by~\eqref{eq:intro_Upsilon}.

Throughout the computation we adopt the convention
that
$\lambda=\pi(\vec{\mathsf{x}})$, $\mu=\overline\pi(\vec{\mathsf{y}})$.

\subsection{Initial Markov kernel}\label{sub:computation_initial_formula}

Our starting point is the formula for
$\Lambda^{N}_{N-1}(\pi(\vec{\mathsf{x}}),\overline\pi(\vec{\mathsf{y}}))$, see
\eqref{eq:traditional_links_Lambda},
under the principal
specialization $x_j=t^{j-1}$, which takes the form (where $\mu\prec\lambda$):
\begin{gather}
	\nonumber
		\Lambda^{N}_{N-1}(\lambda,\mu)
		=
		\psi_{\lambda/\mu}(q,t)\hspace{1pt} x_N^{|\lambda|-|\mu|}\hspace{1pt}
		\frac{
		P_\mu(x_1,\dots,x_{N-1} \mid q,t )
		}
		{P_\lambda(x_1,\dots,x_N \mid q,t)}
		\\
		\nonumber
\hphantom{\Lambda^{N}_{N-1}(\lambda,\mu)}{}
=
		\psi_{\lambda/\mu}
		(q,t)\hspace{1pt}
		t^{(N-1)(|\lambda|-|\mu|)}\hspace{1pt}
		\frac{P_\mu\big(1,t,\dots,t^{N-2}\mid q,t \big)}{P_\lambda\big(1,t,\dots,t^{N-2},t^{N-1}\mid q,t\big)}
		\\
		\nonumber
\hphantom{\Lambda^{N}_{N-1}(\lambda,\mu)}{}
=
		t^{(N-1)(|\lambda|-|\mu|)}
		\prod_{1\le i<j\le N-1}
		\frac{\mathsf{f}\big(q^{\mu_i-\mu_j}t^{j-i}\big)\hspace{1pt}
		\mathsf{f}(q^{\lambda_i-\lambda_{j+1}}t^{j-i})}
		{\mathsf{f}\big(q^{\lambda_i-\mu_j}t^{j-i}\big)\hspace{1pt}
		\mathsf{f}\big(q^{\mu_i-\lambda_{j+1}}t^{j-i}\big)}
		\\
\hphantom{\Lambda^{N}_{N-1}(\lambda,\mu)=}{}
\times
		t^{n(\mu)}\prod_{\square\in \mu}
		\frac{1-q^{a'(\square)}t^{N-1-l'(\square)}}
		{1-q^{a(\square)}t^{l(\square)+1}}
		t^{-n(\lambda)}
		\prod_{\square\in \lambda}
		\frac
		{1-q^{a(\square)}t^{l(\square)+1}}
		{1-q^{a'(\square)}t^{N-l'(\square)}}.		\label{eq:initial_formula}
\end{gather}
Here we used
\eqref{eq:principal_spec} and
\eqref{eq:psi_Pieri}.
Our next steps are devoted to
taking the limit as
$N\to+\infty$
in various parts of the product~\eqref{eq:initial_formula}.

\subsection[Power of t]{Power of $\boldsymbol{t}$}

Let us first consider the overall power of $t$
in \eqref{eq:initial_formula}
which is equal to
$t^{(N-1)(|\lambda|-|\mu|)+n(\mu)-n(\lambda)}$.
Our aim is to express quantities
depending on $\lambda$, $\mu$ through
$\vec{\mathsf{x}}$, $\vec{\mathsf{y}}$.
Adopt the convention
$\mathsf{x}_0=N+m$,
$\mathsf{y}_0=N+m-1$,
$\mathsf{x}_{m+1}=\mathsf{y}_{m+1}=-1$,
so that
\begin{equation*}
	N+m-\mathsf{x}_1-1=\mathsf{x}_0-\mathsf{x}_1-1,\qquad
	N+m-\mathsf{y}_1-2=\mathsf{y}_0-\mathsf{y}_1-1.
\end{equation*}
We have from \eqref{eq:pi_maps_multiplicative_notation}:
\begin{gather*}
	|\lambda|-|\mu|	=
	\sum_{i=0}^{m}
	(N+i)
	(
	( \mathsf{x}_{m-i}-\mathsf{x}_{m-i+1}-1 )-
	( \mathsf{y}_{m-i}-\mathsf{y}_{m-i+1}-1)
)
	\\
\hphantom{|\lambda|-|\mu|}{}
=
	\sum_{i=0}^{m}
	(N+i)
	(
	( \mathsf{x}_{m-i}-\mathsf{x}_{m-i+1} )-
	( \mathsf{y}_{m-i}-\mathsf{y}_{m-i+1} )
)
	\\
\hphantom{|\lambda|-|\mu|}{}
=
	N\sum_{i=0}^{m}
	(
	( \mathsf{x}_{m-i}-\mathsf{x}_{m-i+1} )-
	( \mathsf{y}_{m-i}-\mathsf{y}_{m-i+1} )
	)\\
\hphantom{|\lambda|-|\mu|=}{}
	+
	\sum_{i=0}^{m}i
	(
	( \mathsf{x}_{m-i}-\mathsf{x}_{m-i+1} )-
	( \mathsf{y}_{m-i}-\mathsf{y}_{m-i+1})
	)
	\\
\hphantom{|\lambda|-|\mu|}{}
=
	N+m+
	|\vec {\mathsf{y}}|-|\vec{\mathsf{x}}|.
\end{gather*}
Moreover, from
\eqref{eq:n_lambda_notation} and \eqref{eq:pi_maps_multiplicative_notation}
we have
\begin{gather}
	\nonumber
	n(\mu)-n(\lambda)
	=
	\sum_{j=0}^{m}
	(N+m-j)
	\left[
		\binom{\mathsf{y}_0-\mathsf{y}_{j+1}-(j+1)}2-
		\binom{\mathsf{y}_0-\mathsf{y}_{j}-j}2
	\right]
	\\
\hphantom{n(\mu)-n(\lambda)=}{}
-
	\sum_{j=0}^{m}
	(N+m-j)
	\left[
		\binom{\mathsf{x}_0-\mathsf{x}_{j+1}-(j+1)}2-
		\binom{\mathsf{x}_0-\mathsf{x}_{j}-j}2
	\right]
\nonumber\\
\hphantom{n(\mu)-n(\lambda)}{}
=
	-(N-1)(N+m)
	-\binom{m}{2}
	+|\vec{\mathsf{x}}|\nonumber\\
\hphantom{n(\mu)-n(\lambda)=}{}
+
	\frac12
	\sum_{i=1}^{m}
	(\mathsf{y}_i-\mathsf{x}_i)(\mathsf{x}_i+\mathsf{y}_i+3+2i-2m-2N).
\label{eq:newlabel}
\end{gather}
Indeed, the coefficients
by $\mathsf{x}_i$, $\mathsf{x}_i^2$, $\mathsf{y}_i$, $\mathsf{y}_i^2$, $1\le i\le m$, in the right-hand side of~\eqref{eq:newlabel}
are, respectively,
\begin{equation*}
	N - i + m - \frac{1}{2}
	,\
	-\frac{1}{2}
	,\
	-N + i - m + \frac32
	,\
	\frac{1}{2},
\end{equation*}
which are the same as in the left-hand side.
The free term in the left-hand side is
$-\binom m2+m-N^2-(m-2)N$,
which is readily matched to the free term in the
right-hand side by virtue of
our conventions about
$\mathsf{x}_0$, $\mathsf{y}_0$, $\mathsf{x}_{m+1}$, $\mathsf{y}_{m+1}$.

Let us further simplify the left-hand side of
\eqref{eq:newlabel}.
The $i$-th term in this sum
is rewritten as
\begin{gather*}
	\frac12(\mathsf{y}_i-\mathsf{x}_i)(\mathsf{x}_i+\mathsf{y}_i+3+2i-2m-2N)\\
\qquad{}
	=
	-(N-1)(\mathsf{y}_i-\mathsf{x}_i)
	+
	\frac12(\mathsf{y}_i-\mathsf{x}_i)(-\mathsf{x}_i+\mathsf{y}_i+1)
	+
	(\mathsf{y}_i-\mathsf{x}_i)(\mathsf{x}_i-m+i).
\end{gather*}
Since $\mathsf{y}_i-\mathsf{x}_i$ is equal to $0$ or $-1$,
the quantity
$(\mathsf{y}_i-\mathsf{x}_i)(-\mathsf{x}_i+\mathsf{y}_i+1)$
is identically zero.
Thus, we have
\begin{align*}
	n(\mu)-n(\lambda)=
	-(N-1)(N+m+|\vec{\mathsf{y}}|-|\vec{\mathsf{x}}|)
	+|\vec{\mathsf{x}}|-\binom m2
	+\sum_{i=1}^{m}
	(\mathsf{x}_{i}-m+i)(\mathsf{y}_{i}-\mathsf{x}_{i}).
\end{align*}
We see that
the $N$-dependent terms
in
$(N-1)(|\lambda|-|\mu|)+n(\mu)-n(\lambda)$
cancel out, and the overall
factor containing the
power of $t$
in
$\Lambda^{N}_{N-1}(\pi(\vec{\mathsf{x}}),\overline\pi(\vec{\mathsf{y}}))$
has the form
\begin{equation*}%\label{eq:power_of_t_final}
	t^{-\binom m2+|\vec{\mathsf{x}}|
	+\sum_{i=1}^{m}
	(\mathsf{x}_{i}-m+i)(\mathsf{y}_{i}-\mathsf{x}_{i})}.
\end{equation*}

\subsection{Coarms and colegs}
Addressing the factors in
\eqref{eq:initial_formula}
containing coarms and colegs of $\lambda$ and $\mu$,
we obtain using \eqref{eq:arms_legs}:
\begin{equation*}
	\prod_{\square\in \lambda}\big(1-q^{a'(\square)}t^{N-l'(\square)}\big)
	=
	\prod_{i=1}^{N}
	\prod_{j=1}^{\lambda_i}
	\big(1-q^{j-1}t^{N-i+1}\big)
	=
	\prod_{i=1}^{N}
	\big(t^{N+1-i};q\big)_{\lambda_i}.
\end{equation*}
This product is in the denominator, and a similar
factor
$\prod_{i=1}^{N-1}
\big(t^{N-i};q\big)_{\mu_i}$
appears in the numerator.
Let us show that the contribution coming from
coarms and colegs goes to one, that is,%
\begin{equation}
	\label{eq:coarm_coleg_result}
	\lim_{N\to+\infty}
	\frac{\prod_{i=1}^{N-1}
	\big(t^{N-i};q\big)_{\mu_i}}
	{\prod_{i=1}^{N}
	\big(t^{N+1-i};q\big)_{\lambda_i}}=1.
\end{equation}
To see this, observe that
$\mu_{i}-\lambda_{i+1}$ for all $i=1,\dots,N-1 $ is a nonnegative integer
which does not grow with $N$ as long as $N-i$ is fixed,
see~\eqref{eq:pi_maps_multiplicative_notation}. Therefore,
\begin{gather*}
	\frac{\big(t^{N-i};q\big)_{\mu_i}}{\big(t^{N+1-(i+1)};q\big)_{\lambda_{i+1}}}=
	\frac{\big(t^{N-i};q\big)_{\mu_i}}{\big(t^{N-i};q\big)_{\lambda_{i+1}}}\\
\hphantom{\frac{\big(t^{N-i};q\big)_{\mu_i}}{\big(t^{N+1-(i+1)};q\big)_{\lambda_{i+1}}}}{}
=
	\big(1-t^{N-i}q^{\lambda_{i+1}}\big)\big(1-t^{N-i}q^{\lambda_{i+1}+1}\big)\cdots\big(1-t^{N-i}q^{\mu_i-1}\big) ,
\end{gather*}
where the product in the right-hand side is finite. Since $\lambda_{i+1}$ and $\mu_i$ go to infinity
as $N\to\infty$, this product converges to~$1$.
There is one more factor in the denominator of~\eqref{eq:coarm_coleg_result},
namely,
\begin{equation*}
	\frac{1}{\big(t^N;q\big)_{\lambda_1}}=
	\frac{\big(t^Nq^{N+m};q\big)_\infty}{\big(t^N;q\big)_\infty}.
\end{equation*}
This factor also goes to $1$ as $N\to+\infty$,
and so the limit~\eqref{eq:coarm_coleg_result} is established.

\subsection{Arms and legs}

We now consider the product
in \eqref{eq:initial_formula}
involving arms and legs:
\begin{equation}
	\label{eq:arm_leg_product_initial}
	\prod_{\square\in \mu}
	\frac{1}
	{1-q^{a(\square)}t^{l(\square)+1}}
	\prod_{\square\in \lambda}
	\bigl(1-q^{a(\square)}t^{l(\square)+1}\bigr).
\end{equation}
Recall the notation
$\lambda=\pi(\vec{\mathsf{x}})$,
$\mu=\overline\pi(\vec{\mathsf{y}})$, see~\eqref{eq:pi_maps_multiplicative_notation}.
Observe that the Young diagrams
$\lambda=\pi(\vec{\mathsf{x}})$ and
$\overline\pi(\vec{\mathsf{y}})$ differ only by
adding the first row (cf.\ Figure~\ref{fig:pi_of_x}).
For each box $\square$ in the first row of $\lambda$, the quantity
$l(\square)$ is of order $N$, and so
$1-q^{a(\square)}t^{l(\square)+1}$ is close to~$1$
for large $N$.
Therefore, the product
\eqref{eq:arm_leg_product_initial}
has the same limit as
\begin{equation}	\label{eq:arm_leg_product_proof_1}
	\prod_{\square\in \overline\pi(\vec{\mathsf{y}})}
	\frac{1}
	{1-q^{a(\square)}t^{l(\square)+1}}
	\prod_{\square\in \overline\pi(\vec{\mathsf{x}})}
	\bigl(1-q^{a(\square)}t^{l(\square)+1}\bigr).
\end{equation}

\begin{Proposition}	\label{prop:arm_leg_main_result}
	As $N\to+\infty$, the product \eqref{eq:arm_leg_product_proof_1}
	converges to
	\begin{equation*}%\label{eq:arm_leg_limit_thing}
		\frac{V_{t,q}(\vec{\mathsf{y}})}{V_{t,q}(\vec{\mathsf{x}})}
		\prod_{i=1}^{m}
		\big( 1-q^{m-i}t^{\mathsf{x}_i-m+i}\mathbf{1}_{\mathsf{y}_i=\mathsf{x}_i-1} \big),
	\end{equation*}
	where $V_{t,q}$ is the $(t,q)$-Vandermonde given by \eqref{eq:V_tq}.
\end{Proposition}
\begin{proof}
	Recall the notation
	for the
	$(q,t)$-deformed hook polynomials
	\cite[formulas~(VI.8.1) and (VI.8.1$'$)]{Macdonald1995}:
	\begin{equation*}
		c_\nu(q,t)=
		\prod_{\square\in \nu}\bigl(1-q^{a(\square)}t^{l(\square)+1}\bigr),\qquad
		c'_\nu(q,t)=
		\prod_{\square\in \nu}\bigl(1-q^{a(\square)+1}t^{l(\square)}\bigr),
	\end{equation*}
	and observe that $c_\nu(q,t)=c'_{\nu'}(t,q)$, where $\nu'$ is the
	transposed Young diagram.
	Note that, in multiplicative notation,
	\begin{gather*}
		\overline\pi(\vec{\mathsf{x}})'=
		(N-1)^{N}
		(N-1-\mathsf{x}_m)
		(N-1-\mathsf{x}_{m-1}+1)
		(N-1-\mathsf{x}_{m-2}+2)
		\cdots\\
\hphantom{\overline\pi(\vec{\mathsf{x}})'=}{}
\times (N-1-\mathsf{x}_{1}+m-1)
	\end{gather*}
	and similarly for $\overline\pi(\vec{\mathsf{y}})'$,
	see Figure~\ref{fig:pi_of_x}.
	Adopt the convention
	$\mathsf{x}_{m+j}=\mathsf{y}_{m+j}
	=-j$ for $j=1,2,\dots $.
	Then we may shift the indices $j$ to encode the string $\overline\pi(\vec{\mathsf{x}})'$
	as
	$\overline\pi(\vec{\mathsf{x}})'_j=N+j-1-\mathsf{x}_{m-j}$,
	where $-N\le j\le m-1$, and similarly for $\vec{\mathsf{y}}$.
	The arm and leg lengths do not change under this shifting.
	This shift allows to directly refer to a known identity, the first one in
	\cite[Proposition~3.2]{kaneko1996q}, and write
	\begin{gather*}
		c_{\overline\pi(\vec{\mathsf{x}})}(q,t)
		=
		c'_{\overline\pi(\vec{\mathsf{x}})'}(t,q)
		\\
\hphantom{c_{\overline\pi(\vec{\mathsf{x}})}(q,t)}{}
=
		\frac{(t;t)_{\infty}^{N+m}}
		{
			\prod_{j=-N}^{m-1}
			\big(t^{N+j-\mathsf{x}_{m-j}}q^{m-1-j};t\big)_{\infty}
		}
		\prod_{-N\le i<j\le m-1}
		\frac{\big(q^{j-i}t^{i-j-\mathsf{x}_{m-i}+\mathsf{x}_{m-j}+1};t\big)_{\infty}}
		{\big(q^{j-i-1}t^{i-j-\mathsf{x}_{m-i}+\mathsf{x}_{m-j}+1};t\big)_{\infty}}
		\\
\hphantom{c_{\overline\pi(\vec{\mathsf{x}})}(q,t)}{}
		=
		\frac{(t;t)_{\infty}^{N+m}}
		{
			\prod_{j=-N}^{m-1}
			\big(t^{N+j-\mathsf{x}_{m-j}}q^{m-1-j};t\big)_{\infty}
		}
		\prod_{-N\le i<j\le -1}
		\frac{\big(q^{j-i}t;t\big)_{\infty}}
		{\big(q^{j-i-1}t;t\big)_{\infty}}
		\\
\hphantom{c_{\overline\pi(\vec{\mathsf{x}})}(q,t)=}{}
\times
		\prod_{\substack{-N\le i\le -1 \\ 0\le j\le m-1}}
		\frac{\big(q^{j-i}t^{-j+\mathsf{x}_{m-j}+1};t\big)_{\infty}}
		{\big(q^{j-i-1}t^{-j+\mathsf{x}_{m-j}+1};t\big)_{\infty}}
		\underbrace{
		\prod_{1\le i<j\le m}
		\frac{\big(q^{j-i}t^{i-j-\mathsf{x}_{j}+\mathsf{x}_{i}+1};t\big)_{\infty}}
		{\big(q^{j-i-1}t^{i-j-\mathsf{x}_{j}+\mathsf{x}_{i}+1};t\big)_{\infty}}}
		_{1/V_{t,q}(\vec{\mathsf{x}})}.
	\end{gather*}
	Therefore,
	the product \eqref{eq:arm_leg_product_proof_1}
	becomes
	\begin{gather*}
		\frac{c_{\overline\pi(\vec{\mathsf{x}})}(q,t)}
		{c_{\overline\pi(\vec{\mathsf{y}})}(q,t)}
		=
		\frac{V_{t,q}(\vec{\mathsf{y}})}
		{V_{t,q}(\vec{\mathsf{x}})}
			\prod_{j=1}^{m}
		\frac{
			\big(t^{N+m-j-\mathsf{y}_{j}}q^{j-1};t\big)_{\infty}
		}
		{
			\big(t^{N+m-j-\mathsf{x}_{j}}q^{j-1};t\big)_{\infty}
		}\\
\hphantom{\frac{c_{\overline\pi(\vec{\mathsf{x}})}(q,t)}
		{c_{\overline\pi(\vec{\mathsf{y}})}(q,t)}=}{}
\times
		\prod_{j=1}^{m}
		\prod_{i=-N}^{-1}
		\frac{\big(q^{m-j-i}t^{j-m+\mathsf{x}_{j}+1};t\big)_{\infty}}
		{\big(q^{m-j-i-1}t^{j-m+\mathsf{x}_{j}+1};t\big)_{\infty}}
		\frac
		{\big(q^{m-j-i-1}t^{j-m+\mathsf{y}_{j}+1};t\big)_{\infty}}
		{\big(q^{m-j-i}t^{j-m+\mathsf{y}_{j}+1};t\big)_{\infty}}.
	\end{gather*}
	The first product over $1\le j\le m$
	converges to $1$ as $N\to+\infty$
	thanks to the presence of $t^N$.
	In the second product over $1\le j\le m$, the terms
	where $\mathsf{y}_j=\mathsf{x}_j$ are simply equal to $1$.
	When $\mathsf{y}_j=\mathsf{x}_j-1$, we have
	\begin{gather*}
		\prod_{i=-N}^{-1}
		\frac{\big(q^{m-j-i}t^{j-m+\mathsf{x}_{j}+1};t\big)_{\infty}}
		{\big(q^{m-j-i-1}t^{j-m+\mathsf{x}_{j}+1};t\big)_{\infty}}
		\frac
		{\big(q^{m-j-i-1}t^{j-m+\mathsf{x}_{j}};t\big)_{\infty}}
		{\big(q^{m-j-i}t^{j-m+\mathsf{x}_{j}};t\big)_{\infty}}
	\\
\qquad{}	 =
		\prod_{i=1}^{N}
		\frac{1-q^{m-j+i-1}t^{j-m+\mathsf{x}_{j}}}
		{1-q^{m-j+i}t^{j-m+\mathsf{x}_{j}}}
	 =
		\frac{1-q^{m-j}t^{j-m+\mathsf{x}_{j}}}
		{1-q^{m-j+N}t^{j-m+\mathsf{x}_{j}}},
	\end{gather*}
	and the $N\to+\infty$ limit of this expression is
	$1-q^{m-j}t^{\mathsf{x}_{j}-m+j}$.
	This completes the proof.
\end{proof}

\subsection{Pieri coefficient}

It remains to consider the $N\to+\infty$ limit of the Pieri
coefficient $\psi_{\lambda/\mu}(q,t)$
entering \eqref{eq:initial_formula}.
Recall the convention
$\mathsf{x}_{m+j}=\mathsf{y}_{m+j}
=-j$ for $j=1,2,\dots $,
and encode (shifting the indices)
\begin{equation*}
	\lambda'_i=\pi(\vec{\mathsf{x}})'_i=N+i-\mathsf{x}_{m-i},
	\qquad
	\mu'_j=\overline\pi(\vec{\mathsf{y}})'_j=N-1+j-\mathsf{y}_{m-j}.
\end{equation*}
We can now use
the well-known duality
\cite[proof of formula~(6.24)]{Macdonald1995}
\begin{equation*}
	\psi_{\lambda/\mu}(q,t)=
	\psi'_{\lambda'/\mu'}(t,q),
\end{equation*}
and obtain from~\eqref{eq:psi_prime_Pieri}:
\begin{gather*}
	\psi'_{\lambda'/\mu'}(t,q)
	=
	\prod_{\substack{-N\le i<j\le m-1
	\\
	\lambda'_{i}=\mu'_{i},\, \lambda'_{j}=\mu'_{j}+1}}
	\frac{\big(1-t^{\mu'_{i}-\mu'_{j}} q^{j-i-1}\big)\big(1-t^{\lambda'_{i}-\lambda'_{j}} q^{j-i+1}\big)}
	{\big(1-t^{\mu'_{i}-\mu'_{j}} q^{j-i}\big)\big(1-t^{\lambda'_{i}-\lambda'_{j}} q^{j-i}\big)}
	\\
\hphantom{\psi'_{\lambda'/\mu'}(t,q)=}{}
=
	\prod_{\substack{-N\le i<j\le m-1
	\\
	\mathsf{y}_{m-i}=\mathsf{x}_{m-i}-1,\\
	\mathsf{y}_{m-j}=\mathsf{x}_{m-j}}}
	\frac{\big(1-t^{i-j+\mathsf{y}_{m-j}-\mathsf{y}_{m-i}}
	q^{j-i-1}\big)
	\big(1-t^{i-j+\mathsf{x}_{m-j}-\mathsf{x}_{m-i}} q^{j-i+1}\big)}
	{\big(1-t^{i-j+\mathsf{y}_{m-j}-\mathsf{y}_{m-i}} q^{j-i}\big)
	\big(1-t^{i-j+\mathsf{x}_{m-j}-\mathsf{x}_{m-i}} q^{j-i}\big)}.
\end{gather*}
Split the product over $i<j$ into three parts.
The first part with $-N\le i<j\le -1$
cancels out and is equal to $1$
since
we have $\mathsf{y}_{m-i}=\mathsf{x}_{m-i}$ for all $-N\le i\le -1$.
The second part
with $-N\le i\le -1$ and $0\le j\le m-1$ is also equal to $1$ for the same reason.
The third part with $0\le i<j\le m-1$ equals
\begin{equation}
	\label{eq:psi_prime_in_limit}
	\prod_{\substack{1\le i<j\le m
	\\
	\mathsf{y}_{i}=\mathsf{x}_{i}
	,\,
	\mathsf{y}_{j}=\mathsf{x}_{j}-1}}
	\frac{\big(1-t^{i-j+\mathsf{y}_{i}-\mathsf{y}_{j}}
	q^{j-i-1}\big)
	\big(1-t^{i-j+\mathsf{x}_{i}-\mathsf{x}_{j}} q^{j-i+1}\big)}
	{\big(1-t^{i-j+\mathsf{y}_{i}-\mathsf{y}_{j}} q^{j-i}\big)
	\big(1-t^{i-j+\mathsf{x}_{i}-\mathsf{x}_{j}} q^{j-i}\big)},
\end{equation}
which is independent of $N$. Note that
\eqref{eq:psi_prime_in_limit} is equal to a Pieri coefficient
$\psi'_{(\vec{\mathsf{x}}-\delta_m)/(\vec{\mathsf{y}}-\delta_m)}(t,q)$,
where $\delta_m=(m-1,m-2,\dots,2,1,0)$ is the
staircase partition.

\subsection{Final result}

Putting together all the computations from this section, we see that the
$N\to+\infty$ limit of the Markov kernel $\Lambda^{N}_{N-1}(\pi(\vec{\mathsf{x}}),
\overline\pi(\vec{\mathsf{y}}))$ is the Markov kernel (on $\mathbb{W}_m$)
of the Macdonald noncolliding walks with the absorbing wall at zero:
\begin{gather}
	\nonumber
	\Upsilon_m(\vec{\mathsf{x}},\vec{\mathsf{y}}) =
	t^{-\binom m2+|\vec{\mathsf{x}}|
	+\sum_{i=1}^{m}
	(\mathsf{x}_{i}-m+i)(\mathsf{y}_{i}-\mathsf{x}_{i})}	
	\hspace{1pt}
	\frac{V_{t,q}(\vec{\mathsf{y}})}{V_{t,q}(\vec{\mathsf{x}})}
	\hspace{1pt}
	\psi'_{(\vec{\mathsf{x}}-\delta_m)/(\vec{\mathsf{y}}-\delta_m)}(t,q)
	\\
\hphantom{\Upsilon_m(\vec{\mathsf{x}},\vec{\mathsf{y}})=}{}	
\times
	\prod_{i=1}^{m}
	\big( 1-q^{m-i}t^{\mathsf{x}_i-m+i}\mathbf{1}_{\mathsf{y}_i=\mathsf{x}_i-1}\big)\nonumber
	\\
	\nonumber
\hphantom{\Upsilon_m(\vec{\mathsf{x}},\vec{\mathsf{y}})}{}	
=
	t^{-\binom m2}	
	\hspace{1pt}
	\prod_{1\le i<j\le m}
	\frac
	{
		\big(q^{j-i} t^{\mathsf{x}_i-\mathsf{x}_j-j+i+1};t\big)_{\infty}
		\big( q^{j-i-1}t^{\mathsf{y}_i-\mathsf{y}_j-j+i+1};t \big)_{\infty}
	}
	{
		\big( q^{j-i-1}t^{\mathsf{x}_i-\mathsf{x}_j-j+i+1};t \big)_{\infty}
		\big(q^{j-i} t^{\mathsf{y}_i-\mathsf{y}_j-j+i+1};t\big)_{\infty}
	}
	\\
\hphantom{\Upsilon_m(\vec{\mathsf{x}},\vec{\mathsf{y}})=}{}	\times
\nonumber
	\prod_{\substack{1\le i<j\le m
	\\
	\mathsf{y}_{i}=\mathsf{x}_{i}
	,\,
	\mathsf{y}_{j}=\mathsf{x}_{j}-1}}
	\frac{\big(1-t^{i-j+\mathsf{x}_{i}-\mathsf{x}_{j}+1}
	q^{j-i-1}\big)
	\big(1-t^{i-j+\mathsf{x}_{i}-\mathsf{x}_{j}} q^{j-i+1}\big)}
	{\big(1-t^{i-j+\mathsf{x}_{i}-\mathsf{x}_{j}+1} q^{j-i}\big)
	\big(1-t^{i-j+\mathsf{x}_{i}-\mathsf{x}_{j}} q^{j-i}\big)}
	\\
\hphantom{\Upsilon_m(\vec{\mathsf{x}},\vec{\mathsf{y}})=}{}	\times
	\prod_{i\colon \mathsf{y}_i=\mathsf{x}_i}
	t^{\mathsf{x}_i}
	\prod_{i\colon \mathsf{y}_i=\mathsf{x}_i-1}
	\big(t^{m-i}-q^{m-i}t^{\mathsf{x}_i}\big)
	.\label{eq:Upsilon_final_formula}
\end{gather}
The second expression is obtained by rewriting the
powers of $t$, and expanding the notation of~$V_{t,q}$ and~$\psi'$.
This completes the proof of Theorem~\ref{thm:Macdonald_noncolliding_intro}.

\subsection{Properties of Macdonald noncolliding walks}
\label{sub:absorbing_wall}

Let us show that $\Upsilon_m$ \eqref{eq:Upsilon_final_formula}
indeed possesses an absorbing wall at zero:
\begin{Proposition}
	\label{prop:absorbing_wall}
	Started from any initial configuration,
	the process $\Upsilon_m$ eventually reaches the absorbing state
	$\delta_m = (m-1,m-2,\dots,1,0 )$.
\end{Proposition}
\begin{proof}
	First, observe that if $\mathsf{x}_m=0$, then the term
	$t^{m-i}-q^{m-i}t^{\mathsf{x}_m}$
	(present in
	\eqref{eq:Upsilon_final_formula}
	if \mbox{$\mathsf{y}_m=-1$}) vanishes.
	This means that once the leftmost particle $\mathsf{x}_m$
	reaches $0$, it stays there forever. Moreover, this implies that
	$\delta_m$ is indeed an absorbing state.

	Now, if the process does not eventually reach $\delta_m$,
	then by the monotonicity it will stay an infinite amount of
	time in some configuration $\vec{\mathsf{x}}$ with
	$|\vec{\mathsf{x}}|>|\delta_m|=\binom m2$. However,
	this cannot happen with positive probability
	because $0<t<1$ and
	\begin{equation*}
		\Upsilon_m(\vec{\mathsf{x}},\vec{\mathsf{x}})=
		t^{|\vec{\mathsf{x}}|-\binom m2}<1.
	\end{equation*}
	This completes the proof.
\end{proof}

The walks $\Upsilon_m$ for different $m$ satisfy the following consistency:
\begin{Proposition}
	\label{prop:consistency}
	If $\mathsf{x}_m=\mathsf{y}_m=0$, then
	\begin{equation*}
		\Upsilon_m(\vec{\mathsf{x}},\vec{\mathsf{y}})=
		\Upsilon_{m-1}
		(\vec{\mathsf{x}}\,{}',\vec{\mathsf{y}}\,{}'),
	\end{equation*}
	where $\vec{\mathsf{x}}'\in \mathbb{W}_{m-1}$
	is given by
	$\mathsf{x}_i'=\mathsf{x}_i-1$, $i=1,\dots,m-1 $,
	and similarly for $\vec{\mathsf{y}}\,{}'$.
\end{Proposition}
\begin{proof}
	Follows from the exact formula \eqref{eq:Upsilon_final_formula}
	after necessary simplifications.
\end{proof}

Proposition~\ref{prop:consistency}
allows to view the Markov processes $\Upsilon_m$ for all $m\ge1$
as instances of one and the same Markov process on configurations
of
infinitely many particles on $\mathbb{Z}$.
These configurations must be half-infinite, that is,
there are finitely many particles in $\mathbb{Z}_{\ge0}$
and finitely many holes in $\mathbb{Z}_{<0}$.
If there are $m$ particles away from the packed group
extending to $-\infty$, then the dynamics is governed by
(a suitable shift of) the process $\Upsilon_m$.

\section{Degenerations and limits of Macdonald noncolliding walks}
\label{sec:properties_of_our_walks}

In this section we discuss various
degenerations of the Macdonald noncolliding walks $\Upsilon_m$.

\subsection{Schur degeneration}

When $q=t$, recall that the Macdonald polynomials turn into the Schur polynomials
\eqref{eq:Schur}. This degeneration simplifies our Markov chain, too:

\begin{Proposition}
	\label{prop:Schur_degeneration}
	When $t=q$, the Macdonald noncolliding walks
	$\Upsilon_m$
	\eqref{eq:Upsilon_final_formula}
	on $\mathbb{W}_m$ become
	\begin{gather}\label{eq:Upsilon_Schur}
		\Upsilon_m^{\mathrm{Schur}}(\vec{\mathsf{x}},\vec{\mathsf{y}})
		=
		q^{-\binom m2+(m-1)(|\vec{\mathsf{x}}|-|\vec{\mathsf{y}}|)}
		\!\prod_{1\le i<j\le m}
		\frac
		{
			q^{\mathsf{y}_j}-q^{\mathsf{y}_i}
		}
		{
			q^{\mathsf{x}_j}-q^{\mathsf{x}_i}
		}
		\prod_{i=1}^{m}
		\big( q^{\mathsf{x}_i}\mathbf{1}_{\mathsf{y}_i=\mathsf{x}_i}+
		\big(1-q^{\mathsf{x}_i}\big)\mathbf{1}_{\mathsf{y}_i=\mathsf{x}_i-1}\big)
		.\!\!\!
	\end{gather}
\end{Proposition}
\begin{proof}
	The
	right-hand side of
	\eqref{eq:Upsilon_Schur}
	is a straightforward simplification of \eqref{eq:Upsilon_final_formula}
	at $q=t$.
\end{proof}

This Markov chain may be viewed as the Doob $h$-transform (cf.
\cite{Konig2005, konig2002non})
of $m$ independent random walks on $\mathbb{Z}_{\ge0}$ with transition
probabilities
\begin{equation*}
	\mathrm{Prob}(k\to k)=q^{k},\qquad
	\mathrm{Prob}(k\to k-1)=1-q^{k},\qquad k\in \mathbb{Z}_{\ge0},
\end{equation*}
and absorbing wall at zero.
The corresponding positive harmonic function for the
Markov transition matrix of $m$
independent walks is
\begin{equation*}
	h(\vec{\mathsf{x}})=
	q^{-(m-1)|\vec{\mathsf{x}}|}
	\prod_{1\le i<j\le m}\big(q^{\mathsf{x}_i}-q^{\mathsf{x}_j}\big),
	\qquad \vec{\mathsf{x}}\in \mathbb{W}_{m},
\end{equation*}
and it has the eigenvalue $q^{\binom{m}{2}}$.
The statement about the eigenvalue is equivalent to
\begin{equation*}
	\sum_{\vec{\mathsf{y}}\in \mathbb{W}_m}
	\Upsilon_m^{\mathrm{Schur}}(\vec{\mathsf{x}},\vec{\mathsf{y}})=1,
\end{equation*}
which follows as a $t=q$ specialization of Theorem~\ref{thm:Macdonald_noncolliding_intro}.

The noncolliding walks $\Upsilon_m^{\mathrm{Schur}}$
are somewhat similar to the ones studied in
\cite{BG2011non}. Indeed, the latter
are obtained from $m$ independent random walks on the whole line
$\mathbb{Z}$
by means of the Doob $h$-transform with the $q$-Vandermonde
$\prod_{1\le i<j\le m}\big(q^{\mathsf{x}_i}-q^{\mathsf{x}_j}\big)$, and this
harmonic function has
eigenvalue $1$.
The resulting process studied in \cite{BG2011non}
is invariant with respect to space translations
$\vec{\mathsf{x}}\mapsto \vec{\mathsf{x}}+L$
(where $L\in \mathbb{Z}$ is arbitrary).
In contrast, our walks are not translation invariant and
live on $\mathbb{W}_m\subset \mathbb{Z}_{\ge0}^m$.
As time goes to infinity, our process $\Upsilon_m^{\mathrm{Schur}}$ is
eventually
absorbed at $\delta=(m-1,m-2,\dots,1,0 )\in \mathbb{W}_m$.
It also seems that our noncolliding walks $\Upsilon_m^{\mathrm{Schur}}$
do not admit a limit (while keeping $q$ fixed)
to the process of \cite{BG2011non}.

In \cite{BG2011non} it was shown that the space-time distribution
of the noncolliding $q$-dependent random walks on the whole line $\mathbb{Z}$
is a determinantal point process. The determinantal structure
of the process $\Upsilon_m^{\mathrm{Schur}}$ \eqref{eq:Upsilon_Schur}
will be explored in a forthcoming work.

\subsection{Jack limit}

Fix $\alpha>0$.
When $q=t^\alpha$ and $t\to1$, it is known
\cite[Chapter VI.10]{Macdonald1995}
that Macdonald polynomials turn
into Jack symmetric polynomials.
The parameter $\alpha$ is sometimes denoted by $1/\theta$
(in literature on asymptotic representation theory, for example,
\cite{Kerov1998}),
and is related to the parameter~$\beta$ in random matrix theory as
$\alpha=2/\beta$.

We aim to take the Jack limit of the Markov chain $\Upsilon_m$
\eqref{eq:Upsilon_final_formula}.
As the factors
$t^{m-i}-q^{m-i}t^{\mathsf{x}_i}$ in $\Upsilon_m$
tend to $0$ for fixed $\vec{\mathsf{x}}$,
we need to scale $\vec{\mathsf{x}}$ with $t$.
This scaling necessarily moves the process
away from the absorbing wall.
More precisely, we take the limit as
\begin{equation}
	\label{eq:Jack_limit_transition}
	t\to 1,\qquad L\to+\infty,\qquad
	q=t^{\alpha}, \qquad t^L\to b\in(0,1),
\end{equation}
and shift $\vec{\mathsf{x}}$, $\vec{\mathsf{y}}$ by $L$ as
\begin{equation}
	\label{eq:Jack_limit_transition_2}
	\mathsf{x}_i(L)=\mathsf{X}_i+i(\alpha-1)+L,
	\qquad
	\mathsf{y}_i(L)=\mathsf{Y}_i+i(\alpha-1)+L,
\end{equation}
where
\begin{equation}
	\label{eq:shifted_copies_of_Z}
	\mathsf{X}_i, \mathsf{Y}_i\in \mathbb{Z}-i(\alpha-1),
	\quad
	i=1,\dots,m,\qquad
	\mathsf{X}_i-\mathsf{X}_{i+1}\ge \alpha,\qquad
	\mathsf{Y}_i-\mathsf{Y}_{i+1}\ge \alpha.
\end{equation}
The latter inequalities come from the
strict ordering of
$\vec{\mathsf{x}},\vec{\mathsf{y}}\in \mathbb{W}_m$
\eqref{eq:W_m_space}.
We also have $\mathsf{Y}_i=\mathsf{X}_i$ or $\mathsf{Y}_i=\mathsf{X}_i-1$
for all $i=1,\dots,m $.

\begin{Proposition}	\label{prop:Jack_degeneration}
	Under
	\eqref{eq:Jack_limit_transition}--\eqref{eq:Jack_limit_transition_2},
	the transition probabilities
	$\Upsilon_m(\vec{\mathsf{x}}(L),\vec{\mathsf{y}}(L))$
	converge to
	\begin{gather}
		\label{eq:Upsilon_Jack}
		\Upsilon_m^{\mathrm{Jack}}\big(\vec{\mathsf{X}},\vec{\mathsf{Y}}\big)
		\coloneqq
		b^{m-|\vec {\mathsf{X}}|+|\vec {\mathsf{Y}}|}
		(1-b)^{|\vec {\mathsf{X}}|-|\vec {\mathsf{Y}}|}
		\!\prod_{1\le i<j\le m}\!
		\frac{(\mathsf{X}_i - \alpha \hspace{1pt} \mathbf{1}_{\mathsf{Y}_i=\mathsf{X}_i-1})
		-(\mathsf{X}_j-\alpha \hspace{1pt} \mathbf{1}_{\mathsf{Y}_j=\mathsf{X}_j-1})}
		{\mathsf{X}_i-\mathsf{X}_j}.\!\!\!
	\end{gather}
\end{Proposition}
\begin{proof}
	We have from \eqref{eq:Upsilon_final_formula} for any fixed
	$\vec{\mathsf{x}},\vec{\mathsf{y}}\in \mathbb{W}_m$:
	\begin{gather*}
		\Upsilon_m(\vec{\mathsf{x}}+L,\vec{\mathsf{y}}+L)
		 =
		t^{-\binom m2}	
		\hspace{1pt}
		\prod_{1\le i<j\le m}
		\frac
		{
			\big(q^{j-i} t^{\mathsf{x}_i-\mathsf{x}_j-j+i+1};t\big)_{\infty}
			\big( q^{j-i-1}t^{\mathsf{y}_i-\mathsf{y}_j-j+i+1};t \big)_{\infty}
		}
		{
			\big( q^{j-i-1}t^{\mathsf{x}_i-\mathsf{x}_j-j+i+1};t \big)_{\infty}
			\big(q^{j-i} t^{\mathsf{y}_i-\mathsf{y}_j-j+i+1};t\big)_{\infty}
		}
		\\
\hphantom{\Upsilon_m(\vec{\mathsf{x}}+L,\vec{\mathsf{y}}+L)=}{}
\times
		\prod_{\substack{1\le i<j\le m
		\\
		\mathsf{y}_{i}=\mathsf{x}_{i}
		,\,
		\mathsf{y}_{j}=\mathsf{x}_{j}-1}}
		\frac{\big(1-t^{i-j+\mathsf{x}_{i}-\mathsf{x}_{j}+1}
		q^{j-i-1}\big)
		\big(1-t^{i-j+\mathsf{x}_{i}-\mathsf{x}_{j}} q^{j-i+1}\big)}
		{\big(1-t^{i-j+\mathsf{x}_{i}-\mathsf{x}_{j}+1} q^{j-i}\big)
		\big(1-t^{i-j+\mathsf{x}_{i}-\mathsf{x}_{j}} q^{j-i}\big)}
		\\
\hphantom{\Upsilon_m(\vec{\mathsf{x}}+L,\vec{\mathsf{y}}+L)=}{} \times
		\nonumber
		\prod_{i\colon \mathsf{y}_i=\mathsf{x}_i}
		t^{\mathsf{x}_i+L}
		\prod_{i\colon \mathsf{y}_i=\mathsf{x}_i-1}
		\big(t^{m-i}-q^{m-i}t^{\mathsf{x}_i+L}\big)
		.
	\end{gather*}
	The third line is the only part involving $L$, and
	it turns into
	$b^{m- |\vec {\mathsf{x}}|+|\vec{\mathsf{y}}|}
	(1-b)^{|\vec {\mathsf{x}}| - |\vec{\mathsf{y}}|}$.
	The limits of the other factors are obtained in a standard way,
	for example, see
	\cite[Chapter~VI.10]{Macdonald1995} (and especially formula (VI.10.3)
	and its proof).
	Note that these remaining factors do not depend on the shift $L$, and
	we may thus assume that
	$\mathsf{x}_i=\mathsf{X}_i+i(\alpha-1)$
	and
	$\mathsf{y}_i=\mathsf{Y}_i+i(\alpha-1)$.
	With this notation, we have
	\begin{gather*}
		\frac
		{
			\big(q^{j-i} t^{\mathsf{x}_i-\mathsf{x}_j-j+i+1};t\big)_{\infty}
			\big( q^{j-i-1}t^{\mathsf{y}_i-\mathsf{y}_j-j+i+1};t \big)_{\infty}
		}
		{
			\big( q^{j-i-1}t^{\mathsf{x}_i-\mathsf{x}_j-j+i+1};t \big)_{\infty}
			\big(q^{j-i} t^{\mathsf{y}_i-\mathsf{y}_j-j+i+1};t\big)_{\infty}
		}
		\\
\qquad{}
		\to
		\frac
		{
			\Gamma(\alpha(j-i-1)+\mathsf{x}_i-\mathsf{x}_j-j+i+1)
			\Gamma(\alpha(j-i)+\mathsf{y}_i-\mathsf{y}_j-j+i+1)
		}
		{
			\Gamma(\alpha(j-i)+\mathsf{x}_i-\mathsf{x}_j-j+i+1)
			\Gamma(\alpha(j-i-1)+\mathsf{y}_i-\mathsf{y}_j-j+i+1)
		}
		\\
\qquad\qquad{}
		=
		\frac
		{
			\Gamma(\mathsf{X}_i-\mathsf{X}_j+1-\alpha)
			\Gamma(\mathsf{Y}_i-\mathsf{Y}_j+1)
		}
		{
			\Gamma(\mathsf{X}_i-\mathsf{X}_j+1)
			\Gamma(\mathsf{Y}_i-\mathsf{Y}_j+1-\alpha)
		},
	\end{gather*}
	and
	\begin{gather*}
		\frac{\big(1-t^{i-j+\mathsf{x}_{i}-\mathsf{x}_{j}+1}
		q^{j-i-1}\big)
		\big(1-t^{i-j+\mathsf{x}_{i}-\mathsf{x}_{j}} q^{j-i+1}\big)}
		{\big(1-t^{i-j+\mathsf{x}_{i}-\mathsf{x}_{j}+1} q^{j-i}\big)
		\big(1-t^{i-j+\mathsf{x}_{i}-\mathsf{x}_{j}} q^{j-i}\big)}
		\\
\qquad{}
		\to
		\frac{(\alpha(j-i-1)+i-j+\mathsf{x}_{i}-\mathsf{x}_{j}+1)
		(\alpha(j-i+1)+i-j+\mathsf{x}_{i}-\mathsf{x}_{j})}
		{(\alpha(j-i)+i-j+\mathsf{x}_{i}-\mathsf{x}_{j}+1)
		(\alpha(j-i)+ i-j+\mathsf{x}_{i}-\mathsf{x}_{j})}
		\\
\qquad\qquad{}=
		\frac{
			(\mathsf{X}_i-\mathsf{X}_j+1-\alpha)
			(\mathsf{X}_i-\mathsf{X}_j+\alpha)
		}
		{
			(\mathsf{X}_i-\mathsf{X}_j+1)
			(\mathsf{X}_i-\mathsf{X}_j)
		}.
	\end{gather*}
	One can check by considering four cases
	depending on
	$\mathsf{Y}_i-\mathsf{X}_i\in \left\{ -1,0 \right\}$,
	$\mathsf{Y}_j-\mathsf{X}_j\in \left\{ -1,0 \right\}$
	that the ratio of the gamma functions simplifies as
	\begin{gather*}
		\prod_{1\le i<j\le m}
		\frac
		{
			\Gamma(\mathsf{X}_i-\mathsf{X}_j+1-\alpha)
			\Gamma(\mathsf{Y}_i-\mathsf{Y}_j+1)
		}
		{
			\Gamma(\mathsf{X}_i-\mathsf{X}_j+1)
			\Gamma(\mathsf{Y}_i-\mathsf{Y}_j+1-\alpha)
		}
		\prod_{\substack{1\le i<j\le m
		\\
		\mathsf{Y}_{i}=\mathsf{X}_{i}
		,\,
		\mathsf{Y}_{j}=\mathsf{X}_{j}-1}}
		\frac{
			(\mathsf{X}_i-\mathsf{X}_j+1-\alpha)
			(\mathsf{X}_i-\mathsf{X}_j+\alpha)
		}
		{
			(\mathsf{X}_i-\mathsf{X}_j+1)
			(\mathsf{X}_i-\mathsf{X}_j)
		}
		\\
\qquad{} =
		\prod_{1\le i<j\le m}
		\frac{(\mathsf{X}_i - \alpha \hspace{1pt} \mathbf{1}_{\mathsf{Y}_i=\mathsf{X}_i-1})
		-(\mathsf{X}_j-\alpha \hspace{1pt} \mathbf{1}_{\mathsf{Y}_j=\mathsf{X}_j-1})}
		{\mathsf{X}_i-\mathsf{X}_j},
	\end{gather*}
	which completes the proof.
\end{proof}

The $m$-particle
Markov chain
$\Upsilon_m^{\mathrm{Jack}}$~\eqref{eq:Upsilon_Jack},
where each particle lives on its own shifted copy of~$\mathbb{Z}$
(see~\eqref{eq:shifted_copies_of_Z}),
is a discrete time ``Bernoulli'' analogue
of the $\beta$-nonintersecting
Poisson random walks considered
in~\cite{huang2021beta}.
Indeed, sending $b\to1$, rescaling time from discrete to continuous,
and reversing the direction of jumps from left to right
turns
$\Upsilon_m^{\mathrm{Jack}}$
into the $\beta$-nonintersecting Poisson walks.

When $\alpha=1$ (so the random matrix parameter is $\beta=2$),
the process $\Upsilon_m^{\mathrm{Jack}}$
turns into the process of noncolliding Bernoulli walks
conditioned to never collide.
The trajectory of this Markov process started from an arbitrary
fixed initial configuration is a determinantal point process.
This structure was utilized
in
\cite{GorinPetrov2016universality}
to establish local universality.

\subsection{Hall--Littlewood degeneration
and a continuous time limit}

Let us now take $q=0$. Under this degeneration,
the Macdonald polynomials become
the Hall--Littlewood polynomials
\cite[Chapter~III]{Macdonald1995}.

\begin{Proposition}
	\label{prop:HL_degeneration}
	When $q=0$,
	the Macdonald noncolliding walks $\Upsilon_m$
	\eqref{eq:Upsilon_final_formula} on $\mathbb{W}_m$
	become
	\begin{gather}
			\Upsilon_m^{\mathrm{HL}}
			(\vec{\mathsf{x}},
			\vec{\mathsf{y}}) =
			t^{-\binom m2}	
			\big(t^{\mathsf{x}_m}\mathbf{1}_{\mathsf{y}_m=\mathsf{x}_m}
			+\big(1-t^{\mathsf{x}_m}\big)\mathbf{1}_{\mathsf{y}_m=
			\mathsf{x}_m-1}\big)
			\nonumber\\
\hphantom{\Upsilon_m^{\mathrm{HL}} (\vec{\mathsf{x}},\vec{\mathsf{y}}) =}{}
\times
			\prod_{i=1}^{m-1}
			\big(t^{\mathsf{x}_i}\mathbf{1}_{\mathsf{y}_i=\mathsf{x}_i}
			+t^{m-i}
			\mathbf{1}_{\mathsf{y}_i=\mathsf{x}_i-1}\big)
			\big(1-t^{\mathsf{x}_i-\mathsf{x}_{i+1}-1}
			\mathbf{1}_{\mathsf{y}_i=\mathsf{x}_i-1}
			\mathbf{1}_{\mathsf{y}_{i+1}=\mathsf{x}_{i+1}}\big).\label{eq:Upsilon_HL}
		\end{gather}
\end{Proposition}
\begin{proof}
	A straightforward simplification of \eqref{eq:Upsilon_final_formula}
	at $q=0$.
\end{proof}

In
\eqref{eq:Upsilon_HL},
let us
send $t=(1-\varepsilon)\nearrow 1$ and
scale discrete time
by $\varepsilon$
to continuous time.
This would amount to a Poisson-type
limit transition in our Markov chain
$\Upsilon_m^{\mathrm{HL}}$.

We have
\begin{equation}
	\label{eq:total_jump_rate_HL_cont}
	\Upsilon_m^{\mathrm{HL}}
	(\vec{\mathsf{x}},\vec{\mathsf{x}})=
	t^{-\binom m2+|\vec{\mathsf{x}}|}
	\sim
	1-\big( |\vec{\mathsf{x}}| - \tbinom{m}{2} \big)\varepsilon+
	O\big(\varepsilon^2\big).
\end{equation}
This implies that as $\varepsilon\to0$,
a single step of the Markov chain $\Upsilon_m^{\mathrm{HL}}$
typically does not change~$\vec{\mathsf{x}}$.
However, occasionally, with probability
proportional to $\varepsilon$, a change in $\vec{\mathsf{x}}$
may occur.
All probabilities of order $O\big(\varepsilon^2\big)$
vanish in the scaling limit.

Therefore, a jump in continuous time
can happen in the presence of only one
factor in~\eqref{eq:Upsilon_HL}
proportional to $\big(1-t^k\big)$ for some $k>0$.
Such a factor is associated to a particle
$\mathsf{x}_i$ which has jumped to the left by $1$
while the particle $\mathsf{x}_{i+1}$ has stayed
(if $i=m$, the latter condition
is replaced by $\mathsf{x}_m>0$). This
leads to the jump rate $\mathsf{x}_i-\mathsf{x}_{i+1}-1$,
$i=1,\dots,m $,
where, by agreement, $\mathsf{x}_{m+1}=-1$.
Moreover,
if one particle $\mathsf{x}_i$, $i=1,\dots,m$,
jumps to the left by $1$,
then
any block of particles
$\mathsf{x}_{i-1},\mathsf{x}_{i-2},\dots , \mathsf{x}_{i-r}$
with adjacent indices
can also jump to the left by $1$, at the same rate $\mathsf{x}_i-\mathsf{x}_{i+1}-1$.
Indeed, this is because
any such transition would include the
same factor $1-t^{\mathsf{x}_i-\mathsf{x}_{i+1}-1}\sim (\mathsf{x}_i-\mathsf{x}_{i+1}-1)\hspace{1pt}\varepsilon$.
We can combine these jump events and assign to them the total rate
$i(\mathsf{x}_i-\mathsf{x}_{i+1}-1)$.
When the jump of $\mathsf{x}_i$ happens (at this rate),
then we can additionally select the size of the adjacent block
uniformly at random.
This leads to the following definition of a
continuous time process.

\begin{Definition}
	\label{def:cont_HL}
	Let $\Upsilon_m^{\mathrm{cont}}$ be a continuous
	time Markov process on $\mathbb{W}_m$
	\eqref{eq:W_m_space}
	with jump rates defined
	as follows.
	Attach to each particle $\mathsf{x}_i\in \mathbb{Z}_{\ge0}$,
	$i=1,\dots,m $,
	an independent exponential clock
	of rate $i\hspace{1pt}(\mathsf{x}_i-\mathsf{x}_{i+1}-1)$,
	where, by agreement, $\mathsf{x}_{m+1}=-1$.
	When the clock of $\mathsf{x}_{i}$ rings,
	we additionally
	select an index $j \in \left\{ 1,\dots,i \right\}$
	uniformly at random,
	and all the particles
	$\mathsf{x}_i,\mathsf{x}_{i-1},\dots,\mathsf{x}_j $
	simultaneously
	jump to the left by $1$.
\end{Definition}

The total jump rate from $\vec{\mathsf{x}}$
under the process from
Definition~\ref{def:cont_HL}
is equal to
\begin{equation*}
	m\mathsf{x}_m+(m-1)(\mathsf{x}_{m-1}-\mathsf{x}_m-1)+\dots+
	(\mathsf{x}_1-\mathsf{x}_2-1)=|\vec{\mathsf{x}}| - \tbinom{m}2,
\end{equation*}
which agrees with \eqref{eq:total_jump_rate_HL_cont}.

Therefore, we have established the following
Poisson-type limit transition:
\begin{Proposition}
	\label{prop:HL_t1}
	Let the Hall--Littlewood parameter $t$ be scaled
	as $t=1-\varepsilon$, where $\varepsilon \searrow 0$.
	Let us scale the discrete time as
	$\big\lfloor \varepsilon^{-1}\tau \big\rfloor $, where $\tau\in \mathbb{R}_{\ge0}$
	is the new continuous time parameter.
	Under this scaling, the
	Hall--Littlewood noncolliding walks
	$\Upsilon_m^{\mathrm{HL}}$~\eqref{eq:Upsilon_HL}
	on $\mathbb{W}_{m}$
	turn into the continuous time Markov
	process $\Upsilon_m^{\mathrm{cont}}$
	on $\mathbb{W}_m$ with jump rates given
	by Definition~{\rm \ref{def:cont_HL}} above.
\end{Proposition}

The dynamics $\Upsilon_m^{\mathrm{cont}}$
is somewhat similar to the backwards, inhomogeneous
version of the Hammersley process introduced
in \cite{PetrovSaenz2019backTASEP} in that the
jump rate attached to each particle $\mathsf{x}_i$ is equal to $i$
times the size of the gap behind $\mathsf{x}_i$.
However, the jumping mechanism in
$\Upsilon_m^{\mathrm{cont}}$
is very different from the one in the backwards Hammersley process.

Another observation about
$\Upsilon_m^{\mathrm{cont}}$ is that
it is not clear how to define the ``bulk'' version of the dynamics
living on the full line $\mathbb{Z}$ and preserving a translation invariant
probability distribution on $\left\{ 0,1 \right\}^{\mathbb{Z}}$.
Indeed, the uniformly random
selection of the number of particles which simultaneously
jump at each clock ring is not readily extendable to infinitely many particles
on the full line.
This presents an obstacle to hydrodynamic analysis
of $\Upsilon_m^{\mathrm{cont}}$, even at a heuristic level.

\begin{figure}[th!]
	\centering
	\includegraphics[height=.65\textheight]{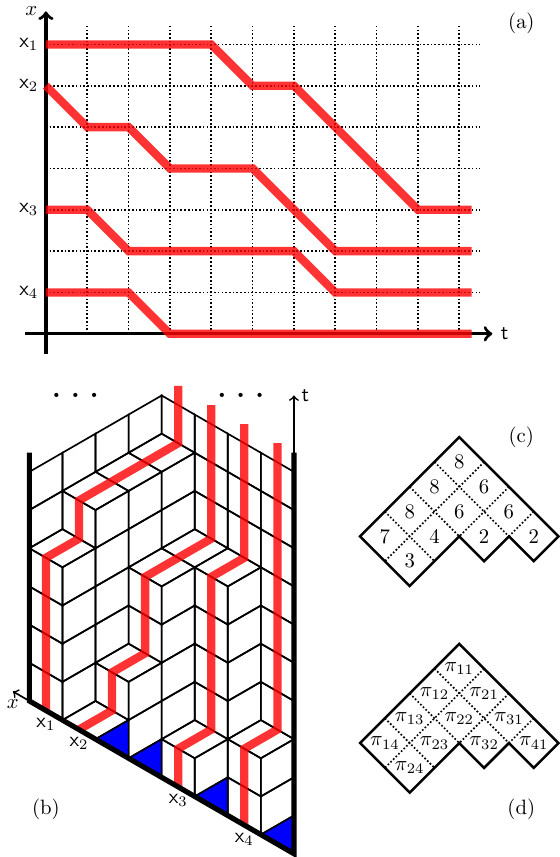}
	\caption{(a) A trajectory of the process
		$\Upsilon_m$ with $m=4$
		and initial configuration
		$\vec{\mathsf{x}}=(7,6,3,1)$;
		(b) A bijective interpretation of the
		trajectory as a lozenge tiling (the shaded
		triangles correspond to $\mathbb{Z}_{\ge0}\setminus
		\vec{\mathsf{x}}$ and
		are removed from the
		strip);
		(c) The corresponding plane partition
		of shape $(4,4,2,1)=(7,6,3,1)-(3,2,1,0)$;
		(d) Encoding the elements of the plane partition which
		must satisfy
		$\pi_{41}\ge0$, $\pi_{32}\ge1$,
		$\pi_{24}\ge 3$, and $\pi_{14}\ge 3.$
	}
	\label{fig:plane_partition}
\end{figure}

\section{Plane partitions}
\label{sec:plane_partitions}

Here we present an interpretation of trajectories of
our noncolliding walks as a certain ensemble of plane partitions
with an arbitrary cascade-like front wall.

\subsection{Bijection to lozenge tilings and plane partitions}

Let $\mathsf{t}\in \mathbb{Z}_{\ge0}$ denote the discrete time
in $\Upsilon_m$, and
let
$\vec{\mathsf{x}}(0)=\vec{\mathsf{x}}
\in \mathbb{W}_m$
be the initial configuration of the process.
Recall from Proposition~\ref{prop:absorbing_wall}
that the final configuration of the process $\Upsilon_m$
is $\delta_m=(m-1,m-2,\dots,1,0 )$.
Via a suitable
affine transform,
let us bijectively map
the trajectory of $\Upsilon_m$
into a lozenge tiling of the vertical strip
of width $\mathsf{x}_1+1$,
see
Figure~\ref{fig:plane_partition}\,(a) and~(b).
The bottom boundary of the vertical strip is encoded by
$\vec{\mathsf{x}}$ in the following way.
Viewing
$\vec{\mathsf{x}}$ as a~particle configuration in
$\mathbb{Z}_{\ge0}$,
each particle $\mathsf{x}_i$
corresponds to a~straight piece in the boundary of slope
$\big({-}1/\sqrt3\big)$,
and each hole
corresponds to cutting a small triangle
out of the strip.\looseness=-1

Due to the absorption (Proposition~\ref{prop:absorbing_wall}),
the lozenge tiling is ``frozen'' far at the top,
with $\mathsf{x}_1+1-m$ tiles of one type on the left followed
by~$m$ tiles of the other type.
Thus, the lozenge tiling contains only finitely many
horizontal lozenges. Therefore, we may view the tiling
as a~graph of a function~$\pi_{i,j}$
defined on cells of a~Young diagram
of size
\begin{equation*}
	\lambda=\vec{\mathsf{x}}-\delta_m=
	(\mathsf{x}_1-m+1,\mathsf{x}_2-m+2,
	\dots,
	\mathsf{x}_{m-1}-1,\mathsf{x}_m ),
\end{equation*}
see Figure~\ref{fig:plane_partition}\,(c) and~(d).
That is, in $\pi_{i,j}$ we have $1\le i\le m$ and
$1\le j\le \lambda_i$.
Since the lozenge tiling cannot have holes inside, this function
must satisfy $\pi_{i,j}\ge\pi_{i+1,j}$ and
$\pi_{i,j}\ge \pi_{i,j+1}$.
Such functions are often called plane partitions, e.g., see~\cite[Chapter~7]{Stanley1999}.

Moreover, due to the sloped bottom boundary
of the lozenge tiling in Figure~\ref{fig:plane_partition}\,(b),
the plane partition must also satisfy
\begin{equation*}%\label{eq:front_wall_condition}
	\pi_{i,\lambda_i}\ge \lambda_i-\mathsf{x}_m
	=
	(\mathsf{x}_i+i)
	-(\mathsf{x}_m+m)
	\qquad \text{for all $1\le i\le m$}.
\end{equation*}
This last condition means that the plane partition
(or rather the corresponding lozenge tiling)
has a cascade-like front wall.
This front wall is encoded by the $\Upsilon_m$'s
initial configuration
$\vec{\mathsf{x}}$ which may be an arbitrary
element of the space $\mathbb{W}_m$
\eqref{eq:W_m_space}.

\subsection{Boltzmann factors}
\label{sub:Gibbs}

Fix $m\ge1$.
For a given fixed initial
configuration,
the space of possible trajectories
of the Markov process
$\Upsilon_m$
is countable.
Here we give a different
characterization of the probability measure
on this space
induced by the transition probabilities
$\Upsilon_m$ \eqref{eq:Upsilon_final_formula}.
Namely, we compute the so-called
Boltzmann factors, that is, the ratios of the
probability weights coming from trajectories related by
an elementary transformation. In this way,
the probability of each given trajectory of $\Upsilon_m$
is proportional to the product of the Boltzmann factors
associated with this trajectory.
Note that such a product
does not depend on the
order of elementary transformations since the
result is proportional to the probability
(the Gibbs weight)
of a trajectory
which we started with.

We need some notation. Fix time $\mathsf{t}\in \mathbb{Z}_{\ge1}$,
and let
$\vec{\mathsf{x}}=\vec{\mathsf{x}}(\mathsf{t}-1)$,
$\vec{\mathsf{y}}=\vec{\mathsf{x}}(\mathsf{t})$,
$\vec{\mathsf{z}}=\vec{\mathsf{x}}(\mathsf{t}+1)$
be three consecutive states of our Markov chain $\Upsilon_m$.
Let us also change $\vec{\mathsf{y}}$ in an elementary
way to $\vec{\mathsf{w}}$ such that for some fixed $1\le k\le m$:
\begin{equation*}
	\mathsf{w}_i=\begin{cases}
		\mathsf{y}_i,&i\ne k,\\
		\mathsf{y}_k-1,& i=k.
	\end{cases}
\end{equation*}
In the lozenge tiling interpretation
of Figure~\ref{fig:plane_partition}\,(b),
the piece of the trajectory
$\vec{\mathsf{x}}\to \vec{\mathsf{y}} \to\vec{\mathsf{z}}$
differs from
$\vec{\mathsf{x}}\to \vec{\mathsf{w}} \to\vec{\mathsf{z}}$
by moving a horizontal lozenge down by $1$.
Note that
$\mathsf{y}_k=\mathsf{x}_k$ and $\mathsf{z}_k=\mathsf{w}_k=\mathsf{y}_k-1$.
In the 3-dimensional interpretation of the lozenge
tiling, this means
removing a $1\times 1\times 1$ box from the stack of boxes.
Let us also denote
\begin{equation*}
	\omega(a,b)\coloneqq
	\frac{
		\big(1-t^{a-b+1}q^{b-1}\big)
		\big(1-t^{a-b} q^{b+1}\big)
	}
	{
		\big(1-t^{a-b+1}q^{b}\big)
		\big(1-t^{a-b}q^{b}\big)
	}
\end{equation*}
to shorten some of the formulas below.

\begin{Proposition}\label{prop:Boltzmann_factor}\samepage
	With the above notation, we have
	\begin{gather}
			\frac{
				\Upsilon_m(\vec{\mathsf{x}},\vec{\mathsf{y}})
				\Upsilon_m(\vec{\mathsf{y}},\vec{\mathsf{z}})
			}
			{
				\Upsilon_m(\vec{\mathsf{x}},\vec{\mathsf{w}})
				\Upsilon_m(\vec{\mathsf{w}},\vec{\mathsf{z}})
			}
			 =
			t
			\prod_{i=1}^{k-1}
			\frac{
				\mathbf{1}_{\mathsf{z}_i=\mathsf{y}_i-1} +
				\omega(\mathsf{y}_i-\mathsf{y}_k,
				k-i)\mathbf{1}_{\mathsf{z}_i=\mathsf{y}_i}
			}
			{
				\mathbf{1}_{\mathsf{y}_i=\mathsf{x}_i-1}
				+
				\omega(\mathsf{x}_i-\mathsf{x}_k,k-i)
				\mathbf{1}_{\mathsf{y}_i=\mathsf{x}_i}
			}
			\nonumber\\
\hphantom{\frac{
				\Upsilon_m(\vec{\mathsf{x}},\vec{\mathsf{y}})
				\Upsilon_m(\vec{\mathsf{y}},\vec{\mathsf{z}})
			}
			{
				\Upsilon_m(\vec{\mathsf{x}},\vec{\mathsf{w}})
				\Upsilon_m(\vec{\mathsf{w}},\vec{\mathsf{z}})
			}
			 =}{}
			\times
			\prod_{j=k+1}^{m}
			\frac{
				\mathbf{1}_{\mathsf{y}_j=\mathsf{x}_j}
				+
				\omega(\mathsf{x}_k-\mathsf{x}_j,j-k)
				\mathbf{1}_{\mathsf{y}_j=\mathsf{x}_j-1}
			}
			{
				\mathbf{1}_{\mathsf{z}_j=\mathsf{y}_j}
				+
				\omega(\mathsf{y}_k-\mathsf{y}_j-1,j-k)
				\mathbf{1}_{\mathsf{z}_j=\mathsf{y}_j-1}
			}.\label{eq:Boltzmann_factor}
		\end{gather}
\end{Proposition}
\begin{proof}
	We use the formula 	
	\begin{gather*}
		\Upsilon_m(\vec{\mathsf{x}},\vec{\mathsf{y}})=
		t^{-\binom m2+|\vec{\mathsf{x}}|
		+\sum_{i=1}^{m}
		(\mathsf{x}_{i}-m+i)(\mathsf{y}_{i}-\mathsf{x}_{i})}	
	\\
\hphantom{\Upsilon_m(\vec{\mathsf{x}},\vec{\mathsf{y}})=}{}\times
		\frac{V_{t,q}(\vec{\mathsf{y}})}{V_{t,q}(\vec{\mathsf{x}})}
		\hspace{1pt}
		\psi'_{(\vec{\mathsf{x}}-\delta_m)/(\vec{\mathsf{y}}-\delta_m)}(t,q)
		\prod_{i=1}^{m}
		\big( 1-q^{m-i}t^{\mathsf{x}_i-m+i}\mathbf{1}_{\mathsf{y}_i=\mathsf{x}_i-1} \big)
	\end{gather*}
	for the transition probability.
	Clearly, in their combination
	in the left-hand side of~\eqref{eq:Boltzmann_factor}
	all factors $V_{t,q}$ cancel out.
	Next, recall that $|\vec{\mathsf{y}}|=|\vec{\mathsf{w}}|+1$,
	and this gives rise to the factor $t$ in front of
	the right-hand side of \eqref{eq:Boltzmann_factor}.
	Next, one can readily check that all the
	factors coming from
	\begin{equation*}
		t^{-\binom m2
		+\sum_{i=1}^{m}
		(\mathsf{x}_{i}-m+i)(\mathsf{y}_{i}-\mathsf{x}_{i})}
		\prod_{i=1}^{m}
		\big( 1-q^{m-i}t^{\mathsf{x}_i-m+i}
		\mathbf{1}_{\mathsf{y}_i=\mathsf{x}_i-1} \big)
	\end{equation*}
	in the left-hand side of \eqref{eq:Boltzmann_factor}
	cancel out, too.
	Finally, we are left with
	the ratio
	\begin{align*}
		\frac{
		\psi'_{(\vec{\mathsf{x}}-\delta_m)/(\vec{\mathsf{y}}-\delta_m)}(t,q)
		\psi'_{(\vec{\mathsf{y}}-\delta_m)/(\vec{\mathsf{z}}-\delta_m)}(t,q)
		}
		{
		\psi'_{(\vec{\mathsf{x}}-\delta_m)/(\vec{\mathsf{w}}-\delta_m)}(t,q)
		\psi'_{(\vec{\mathsf{w}}-\delta_m)/(\vec{\mathsf{z}}-\delta_m)}(t,q)}.
	\end{align*}
	Recalling the definition of
	$\psi'$ \eqref{eq:psi_prime_Pieri}, we may rewrite this ratio as
	\begin{gather*}
		\prod_{\substack{1\le i<j\le m
		\\
		\mathsf{y}_{i}=\mathsf{x}_{i}
		,\,
		\mathsf{y}_{j}=\mathsf{x}_{j}-1}}
		\omega(\mathsf{x}_i-\mathsf{x}_j,j-i)
		\prod_{\substack{1\le i<j\le m
		\\
		\mathsf{z}_{i}=\mathsf{y}_{i}
		,\,
		\mathsf{z}_{j}=\mathsf{y}_{j}-1}}
		\omega(\mathsf{y}_i-\mathsf{y}_j,j-i)
		\\\qquad{}
		\times
		\prod_{\substack{1\le i<j\le m
		\\
		\mathsf{w}_{i}=\mathsf{x}_{i}
		,\,
		\mathsf{w}_{j}=\mathsf{x}_{j}-1}}
		\frac{1}{\omega(\mathsf{x}_i-\mathsf{x}_j,j-i)}
		\prod_{\substack{1\le i<j\le m
		\\
		\mathsf{z}_{i}=\mathsf{w}_{i}
		,\,
		\mathsf{z}_{j}=\mathsf{w}_{j}-1}}
		\frac{1}{\omega(\mathsf{w}_i-\mathsf{w_j},j-i)}.
	\end{gather*}
	Here all terms where neither $i$ nor $j$ is equal to $k$
	cancel out. When $j=k$, we may only get nontrivial
	contributions from the second and the third products,
	and when $i=k$, nontrivial contributions may only come
	from the first and the fourth products. In the
	fourth product, we use
	$\mathsf{w}_k-\mathsf{w}_j=\mathsf{y}_k-\mathsf{y}_j-1$.
	This completes the proof.
\end{proof}

In the Schur case $t=q$,
we have $\omega(a,b)=1$,
so the expression~\eqref{eq:Boltzmann_factor}
reduces simply to~$q$.
In this way, adding a $1\times 1\times 1$
box multiplies the probability weight of a
lozenge tiling by~$q$, so the
whole probability of a tiling is
proportional to $q^{\mathsf{vol}}$,
where $\mathsf{vol}$ is the volume under the
corresponding 3-dimensional surface. Note that
this volume is
the same as
the sum of the entries of the plane partition
as in Figure~\ref{fig:plane_partition}\,(c) and~(d).

Gibbs measures
on lozenge tilings of various infinite regions
in which the probability weight of a tiling is
proportional to $q^{\mathsf{vol}}$
have been widely studied, see, for example,
\cite{BMRT2009BackWall,okounkov2003correlation,Okounkov2005}.
Most well-known ensembles of such lozenge tilings
are solvable by means of Schur processes, and feature
an arbitrary back wall. Our ensemble of~$q^{\mathsf{vol}}$ weighted
lozenge tilings possesses a different kind of boundary conditions,
namely, an arbitrary cascade front wall, as
seen in Figure~\ref{fig:plane_partition}\,(b).

\subsection*{Acknowledgments}

I am grateful to Alexei Borodin,
Grigori Olshanski, and Mikhail Tikhonov for fruitful discussions,
and to the anonymous referees for
helpful remarks.
The work was partially supported by the NSF
grant DMS-1664617, and the
Simons Collaboration Grant for Mathematicians 709055.
This material is based upon work supported by the National
Science Foundation under Grant No. DMS-1928930 while LP
participated in the program
``Universality and Integrability in random matrix theory and Interacting Particle Systems''
hosted by the Mathematical
Sciences Research institute in Berkeley, California, during
the Fall 2021 semester.

\pdfbookmark[1]{References}{ref}
\LastPageEnding

\end{document}